\documentclass[10pt]{amsart}
\usepackage{amsmath,amssymb,amscd,latexsym,eufrak,amsfonts,amsthm,mathrsfs}
\input xy
\xyoption{all} 

\usepackage{array,latexsym,color,a4wide,float,enumerate}
\usepackage{graphicx,epsfig}
\usepackage[all]{xy}
\usepackage{color,mathrsfs}
\usepackage{longtable}

\newcommand{\ol}{\overline}
\newcommand{\kn}{K_{\ov N}}
\newcommand{\ble}{\begin{lemma}}
\newcommand{\ele}{\end{lemma}}
\newcommand{\ti}{\tilde}

\theoremstyle{plain} 
\newtheorem{theorem}{Theorem}[section]
\newtheorem{lemma}[theorem]{Lemma}
\newtheorem{proposition}[theorem]{Proposition}
\newtheorem{corrolary}[theorem]{Corollary}

\theoremstyle{definition}
\newtheorem{definition}[theorem]{Definition}
\newtheorem{example}[theorem]{Example}
\newtheorem{noname}[theorem]{}
\newtheorem{remark}[theorem]{Remark}
\newtheorem{construction}[theorem]{Construction}

\theoremstyle{remark}
\newtheorem*{smallremark}{Remark}
\newtheorem*{note}{Note}

\newtheorem{case}{Case} \makeatletter \@addtoreset{case}{theorem}\makeatother
\newtheorem{claim}{Claim} \makeatletter \@addtoreset{claim}{theorem}\makeatother

\newcommand{\bthm}{\begin{theorem}}
\newcommand{\bprop}{\begin{proposition}}
\newcommand{\blem}{\begin{lemma}}
\newcommand{\bcor}{\begin{corrolary}}
\newcommand{\brem}{\begin{remark}}
\newcommand{\bdfn}{\begin{definition}}
\newcommand{\bitem}{\begin{itemize}}
\newcommand{\bex}{\begin{example}}
\newcommand{\bno}{\begin{noname}}
\newcommand{\bsrem}{\begin{smallremark}}
\newcommand{\bnot}{\begin{note}}
\newcommand{\bcon}{\begin{construction}}
\newcommand{\bca}{\begin{case}}
\newcommand{\bcl}{\begin{claim}}

\newcommand{\ecl}{\end{claim}}
\newcommand{\eca}{\end{case}}
\newcommand{\econ}{\end{construction}}
\newcommand{\enot}{\end{note}}
\newcommand{\esrem}{\end{smallremark}}
\newcommand{\eno}{\end{noname}}
\newcommand{\eex}{\end{example}}
\newcommand{\eitem}{\end{itemize}}
\newcommand{\ethm}{\end{theorem}}
\newcommand{\eprop}{\end{proposition}}
\newcommand{\elem}{\end{lemma}}
\newcommand{\ecor}{\end{corrolary}}
\newcommand{\erem}{\end{remark}}
\newcommand{\edfn}{\end{definition}}

\newcommand{\ba}{\begin{array}}
\newcommand{\ea}{\end{array}}
\newcommand{\benum}{\begin{enumerate}}
\newcommand{\eenum}{\end{enumerate}}

\newcommand{\ov}{\overline}

\newcommand{\disjointslash}{\mathinner{\mskip-1mu\raise1pt\hbox{{\small /}}\mskip-13mu\raise-1.58pt\hbox{{\Large \textvisiblespace}}\mskip0mu}}

\def\m1{$(-1)$-curve}
\def\m1s{$(-1)$-curves}

\def\ds{\displaystyle\sum}

\def\8{\infty}

\def\PP{\mathbb{P}}

\def\C{\mathbb{C}}
\def\Z{\mathbb{Z}}

\def\Q{\mathbb{Q}}

\def\ovk{\overline\kappa}
\def\kp{\kappa}
\def\P{{\mathbb P}}
\def\ks{K_{\ov S}}
\def\kn{K_{\ov N}}

\def\:{\colon}

\newcommand{\To}{\rightarrow}

\def\med{\medskip}

\newcommand{\noin}{\noindent}

\def\Bk{\operatorname{Bk}}

\def\Supp{\operatorname{Supp}}
\def\Pic{\operatorname{Pic}}

\begin{document}

\title[Some properties of $\C^*$ in $\C^2$. ]{Some properties of $\C^*$ in $\C^2$. }

\begin{abstract}{ We consider plane curves isomorphic to $\mathbb{C}^*$. We prove that with one exception the branches at infinity can be separated by an automorphism of $\mathbb{C}^2$. We also give a bound for selfintersection number of the resolution curve.}\end{abstract}

\author[M.Koras and P.Russell]{Mariusz Koras and Peter Russell}
\address{M.Koras: Institute of Mathematics, University of Warsaw and Peter Russell: Department of Mathematics, McGill University, Montreal, Canada. }
\thanks{The first author was supported by Polish Grant N N201 608640.
The second author was supported by a grant from NSERC, Canada.}
\email{koras@mimuw.edu.pl, russell@math.mcgill.ca}
\maketitle

\setcounter{section}{-1}
\section{Introduction}

\bno\label{setup} Let $U$ be a closed algebraic curve in
$\C^2$ isomorphic to $\C^*$. Let $\ov{U}$ be the closure of $U$ in
$\mathbb{P}^2$. By $L_{\infty}$ we denote the line at infinity in
$\mathbb{P}^2$.  Let $\Phi\colon \ov{S}'\rightarrow \mathbb{P}^2$
be the resolution of $\ov{U}$. By this we mean that $\Phi^{-1}$ is the minimal
sequence of blow ups such that the reduced inverse image of the
divisor $\ov{U}+L_{\infty}$ is an SNC-divisor. Let $E'$ be the
proper transform of $\ov{U}$ and let
$D'=\Phi^{-1}(L_{\infty})_{red}$. Let $L_{\infty}'$ denotes the
proper transform of the line $L_{\infty}$ in $\ov{S}'$.  Let $\Psi\colon \ov{S}'\rightarrow
\ov{S}$ be the {\it NC-minimalization of the divisor $D'$ with respect to $E'$}, i.e.,
$\Psi$ is the successive contraction of possibly $L_{\infty}'$ and then more $D'$-components such that
$\Psi(D'+E')$ is a SNC-divisor and  each $(-1)$-component of $\Psi(D')$ is
a branching component of $\Psi(D'+E')$. We put $D=\Psi(D')$, $E=\Psi(E')$. Let
$S=\ov{S}\setminus D$.  Of course $S\simeq\mathbb{C}^2$. We note that $E'\cdot D'=E\cdot D=2$. Since $D$ has connected support, $D+E$ is not a chain.
\eno

\med Embeddings of $\C^*$ into $\C^2$ can be divided into two classes. The first class consists of embeddings which admit  {\it a good asymptote}, see \ref{asymptote},  the second class consists of those without any good asymptote. The embeddings from the first class are completely classified in \cite{C-NKR}.

\begin{definition}\label{asymptote} We say that a  rational curve $L\in \PP^2$ is a \textit{good asymptote} of $U$ if
 $L\cap \C^2\simeq \C^1$,  and $L$ meets $\ov U$
 at most once  at finite distance, i.e., $L\cdot U\leq 1$.
 \end{definition}

 Notice that this definition differs slightly from the definition in \cite{C-NKR}, but the two definitions are equivalent up to an isomorphism of $\C^2$.

\med The main results of this article are Corollary 2.5 and Theorem 4.16. Corollary 2.5 gives a bound for $E^2$ and  for $(\ks+D+E)^2$ in the case where $U$ does not admit a good asymptote. Theorem 4.16 says that in the case of no good asymptote the branches of $\ov U$ at infinity can be separated by an automorphism of $\C^2$. It follows from the classification given in \cite{C-NKR} that with one exception this is also true in case where $U$ admits a good asymptote. Hence throughout the paper we assume that $U$ does not have a good asymptote.

\med Another remarkable property  is proved in \cite{Kor}.

\begin{theorem} $\ovk(\ov S\setminus E)=-\infty.$
\end{theorem}

\med Theorem 0.3  and a theorem of Coolidge imply that $U$ can be transformed into a line in $\C^2$ by a birational authomorphism of $\C^2$, see \cite{KM}.

\section{Preliminaries}
In the article we use several notions and results from the
theory of open algebraic surfaces.  We refer
the reader to \cite{M} for any undefined terms here. We will also use some
results from T. Fujita's paper \cite{Fu1}, particularly $\S$ 3.

\begin{noname}\label{1.1}  Let $\ov{M}$ be a complete, non-singular surface and
$T=\sum_{i=1}^n m_iT_i$ a  divisor on $\ov{M}$ with $T_1,\dots , T_n$
distinct, irreducible curves.

 (i) We write $\sim$ for linear equivalence of integral divisors. We write
$\equiv$ for numerical equivalence of divisors, both over $\Z$ and
over $\Q$.

\med (ii) We call $T$ a simple normal crossing  divisor (an $SNC$-divisor) if $T$ is reduced, all its components are smooth
 and at most two of  them
meet at any point, and if so, transversally.

\med (iii) A $(b)$-curve on $\ov{M} $ is a curve $L\simeq \PP^1$ with $L^2=b$.

\med (iv) An $SNC$-divisor $T$ is $NC$-minimal if every (-1)-component of $T$ is a branching component.

\med (v) We call $Q(T)=(T_i\cdot T_j)_{1\leq i,j\leq n }$ the intersection matrix of a reduced $T$ and put $d(T)=det(-Q(T))$. We put $d(T)=1$ if  $\Supp(T)=\emptyset.$

\med (vi) A divisor $R$ is called {\it contractible} if it is the minimal resolution divisor of a quotient singular point. Hence $R$ is a chain composed of smooth rational curves $R_i$ such that $R_i^2\leq-2$  or $R$ is a fork of of smooth rational curves with branches of type $(2,2,n), (2,3,3), (2,3,4)$ or $(2,3,5)$ with and negative branching component.

\end{noname}


\

\bno\label{1.3} For the definition  of {\it twig, tip, bark} of  a divisor   and their properties we refer to
 \cite[\S3]{Fu1} and \cite[\S2, section 3]{M}. We recall only the definition of a {\it capacity} of a rational chain. Let
\begin{displaymath}
R=R_1+\cdots +R_s
\end{displaymath}
be a chain of smooth rational curves with dual graph
\begin{displaymath}
\xymatrix{
*[r]{\otimes}\ar@{}[d]|{b_1}\ar @{-} [r] & {\dots}\ar @{-} [r] & *[l]{\otimes}\ar@{}[d]|{b_s}\\
{} & {} & {} }
\end{displaymath}

 Suppose  that $R$ is admissible, i.e., that $b_i=R^2_i\leq -2,\ i=1,\cdots, s$.  We recall that $Q(R)$ is negative definite and $d(R)\geq 2$. We put $e(R)=\frac{d(R_2+\dots +R_s)}{d(R)}$.
\

\end{noname}

\med\noin\textbf{1.2.1.} {\it If $R_1^2=-k$, then $e(R)\geq\frac{1}{k}$.}


\med We recall
\bno
\med\noin{\it Let $T$ be a connected $NC$-minimal divisor consisting of smooth rational curves. Assume $T$ is not a contractible divisor and let $T_1,\dots, T_s$ be all the maximal twigs of $T$. Then $\Bk(T)^2=-\sum e(T_i)$.}
\eno

 \ble There is no curve $C\subset \ov S$ such that $C\cap S\simeq \C^1$ and $C\cdot E\leq 1$. \ele

    \begin{proof}   Let $L$ be the proper transform of $C$ in $\PP^2$. Then $L$ is a good asymptote of $\ov U$; a contradiction.\end{proof}

    \bcor\label{S-E minimal} If $E^2\neq -1$,  then the pair
  $(\ov S, D+E) $ is  relatively minimal (see \cite[\S2, section 3]{M}) i.e. $\ks+D+E\equiv P+Bk(D+E)$ is the Zariski decomposition, where $P=(\ks+D+E)^+$.
    \ecor

\ble\label{kappa S-E}
\benum
\item[{\rm (i)}]   $\ovk(\ov S\setminus (D+E))\geq 0$.

\item[{\rm (ii)}] $(\ks+D+E)^2<3$.
\eenum
\ele

\begin{proof}   (i) If $\ovk(\ov S\setminus (D+E))=-\infty$, then $|\ks+D+E|=\emptyset$, which implies $E\cdot D\leq 1$, but $E\cdot
D=2$. (ii) By \cite{Mi}, $(\ks +D+E)^2\leq 3(\chi(S\setminus
E)+\frac{1}{12}N^2)$, where $N=(\ks +D+E)^-$ is the negative part
in the Zariski decomposition of the divisor $\ks+D+E$. If $N\neq 0$
we are done since $\chi(S\setminus E)=1$. Suppose that $N=0$. It
follows that the divisor $D+E$ has no twigs (see \cite[\S3]{Fu1},\cite[\S2, section 3]{M}). Hence $D$ is a chain and $E$ meets the tips of $D$.
If $D$ has only one component, then $D^2=1$. Hence in all cases no component of $D$ is a (-1)-curve since $D$ is NC-minimal w.r.t $E$, see 0.1. Clearly $D$ is not an
admissible chain, see 1.2. Therefore there exists a component $D_1$  of $D$ such that $D_1^2\geq 0$. By some blowing up and down within $D$ we can transform $D$ into a chain $\Delta$ with a tip $\Delta_1$ such that $\Delta_1^2=0$ and $E\cdot \Delta_1=1$. The linear system $\vert \Delta_1 \vert$ induces a
$\mathbb{C}$-ruling of $S$ with $E$ as a 1-section. The proper transform in $\PP^2$ of a general member of the system  is a good asymptote of $U$, contrary to our assumption.
\end{proof}


\bno\label{epsilon} Write
$$(\ks+D+E)^2=2-\varepsilon$$

\noin where $\varepsilon\geq 0.$
 We have $\ks\cdot (\ks+D+E)=2-\varepsilon.$ We put $\gamma=-E^2$.
 \eno
 \ble\label{gamma>0} $\gamma>0$.
 \ele
 \begin{proof} Suppose that $\gamma\leq 0$ After blowing up over one of the points in $E\cap D$ we may assume that $E^2=0$. Therefore $U$ is a fiber of a $\C^*$-ruling of $\C^2$. There is a singular fiber with an irreducible component isomorphic to $\C$. This is a good asymptote of $U$ and we reach a contradiction.\end{proof}

\bno
Let $R_ 1,\cdots ,R_s$ be all maximal twigs of $D+E$. Let $e_i=e(R_i)$.
\eno

\ble\label{sum ei}  $-(Bk(D+E))^2=\sum e_i\leq 1+\varepsilon$.
\ele
\begin{proof} Suppose that $\gamma\neq 1$. Then, by \ref{S-E minimal}, the pair $(\ov S,D+E)$ is minimal. Let $\ks+D+E=P+\Bk(D+E)$  be the Zariski decomposition. By Langer's version \cite{L} of the Kobayashi inequality, see \ref{BMY}, $0\leq P^2\leq 3\chi(S\setminus E)=3$. We have $P^2=(\ks+D+E)^2-(\Bk(D+E))^2=2-\varepsilon+\sum e_i$, so we are done. Suppose that $E^2=-1$. We pass to a minimal model of the pair $(\ov S,D+E)$. In view of 1.4, we possibly contract $E$ and further components of $D+E$, but we do not  touch any of the maximal twigs of $D+E$. To the resulting divisor we apply the Kobayashi inequality and get the result.
\end{proof}

\ble\label{1.13} $\varepsilon\leq 3$.
\ele
\begin{proof} We claim that $\ov S$ is not isomorphic to a relatively minimal rational surface. By \ref{gamma>0}, $E^2< 0$.  Hence $\ov S$ is not isomorphic to $\PP^2$, and if it is isomorphic to a Hirzebruch surface, then $E$ is the only negative curve in $\ov S$. Also $D=D_1+D_2$ has two irreducible components since the irreducible components of $D$ generate $\Pic(\ov S)$ freely. Now $D_1^2\geq 0, D_2^2\geq 0$, and since $E\cdot (D_1+D_2)=2$ we may assume $E\cdot D_1 \leq 1$, say. After some blowing up we may assume that $D_1^2=0$ and then the proper transform of a general member of the system $\vert D_1\vert $ is a good asymptote of $U$, in contradiction to our assumption. \\
Suppose that $\varepsilon \geq 4$. We have $(\ks+D)\cdot(\ks+E)=2-\varepsilon -\ks\cdot E+(\ks+D)\cdot E=4-\varepsilon \leq 0$. Since $E\cdot D=2$, $\vert E+\ks+D\vert\neq \emptyset$. By  Fujita's Theorem \cite {Fu2} there exists $m$ such that $\vert E+m(\ks+D)\vert \neq\emptyset$, but $\vert E+(m+1)(\ks+D)\vert =\emptyset$. We write
$$E+m(\ks+D)=\sum A_i$$
with each $A_i$ reduced and irreducible. We have $\vert A_i+D+\ks\vert=\emptyset$ for every $i$. By a standard argument, \cite [2.1, 2.2] {Ru2} for example, $A_i$ is a smooth rational curve and $A_i\cdot D\leq 1$. This implies $A_i\neq E$. Since $\ov S$ is not a relatively minimal surface we may assume that $A_i^2<0$ for every $i$. (If $A_i^2\geq 0$ we replace $A_i$ by a suitable singular member of the linear system $\vert A_i \vert$). We obtain $-2\geq E\cdot(\ks+E)+m(\ks+D)\cdot(\ks+E)=\sum A_i\cdot(\ks+E)$. Hence there exists $A_1$ such that $A_1\cdot (\ks+E)<0$. It follows that $A_1\cdot \ks<0$. Hence $A_1^2=-1$ and $A_1\cdot E=0$. Since $A_1\cdot D\leq 1$, $A_1$ is not a branching component of $D$. Since $A_1\cdot E=0$ it follows from the NC-minimality of $D$ w.r.t. $E$ (see 0.1) that $A_1$ is not a component of $D$. Now the proper transform of $A_1$ in $\PP^2$ is a good asymptote of $U$, a contradiction.

\end{proof}

\ble Let $\ov M $ be a smooth projective surface. Let $r$ be the rank the Neron-Severi group $NS(\ov M)$. Then for any
set $C_1,\cdots, C_r$ of distinct irreducible curves  in $\ov M$ the matrix $[C_i\cdot C_j]$ is not negative definite.
\ele

\begin{proof} Suppose it is. Then in particular $C_1,\dots,C_r$ are independent in $NS(\ov M)\otimes
\mathbb{Q}$ and hence form a basis of $NS(\ov M)\otimes
\mathbb{Q}$. We reach a contradiction with the Hodge Index Theorem.\end{proof}

We will use an  inequality of Bogomolov-Miyaoka-Yau type (simply the BMY-inequality) proved by R. Kobayashi,
S. Nakamura and F. Sakai ,
\cite[Lemma 8 and Corollary 9]{GM} (see also \cite[Chapter 2, Theorem 6.6.2]{Mi}). We state it as follows.

\ble\label{BMY}
Let $\ov{X}$ be a smooth projective surface and let $D$ be an SNC-divisor on $\ov X$. Let $D_1,\dots ,D_k$ be the connected components of $D$ which are contractible divisors, see 1.1(vi). Let $G_i, 1\leq i\leq k$, be the local fundamental group at the singular point obtained by the contraction of $G_i$ to point. Assume that the pair $(X,D)$ is almost minimal.  Suppose that
$\ovk(X\setminus D)\geq 0$. Then

\[
((K_{\ov X}+D)^+)^2\leq 3(\chi(X\setminus D)+\sum_{i=1}^k\frac{1}{|G_i|})  ,
\]

\ele
 The original BMY-inequality was proved in case $\ovk(X\setminus D)=2$. A. Langer \cite{L}
has extended it to the case $\ovk(X\setminus D)=0,1$.

\section{Basic inequality}

\bno Let $\psi\colon \ov{S}\To \ov{N}$ be a 2-reduction of the divisor $D$ with respect to $E$, i.e., $\psi$ is a
sequence of successive contractions of $(-1)$-curves in $D$ meeting $E$ ( and its successive images) once and such that the divisors $T=\psi (D)$ and $E_0=\psi(E)$ satisfy the
following:

(i) $T$ is an NC-divisor.

(ii) for any (-1)-component $T_i$ of $T$, $T_i\cdot E_0\geq 2 $ or $T_i$ is a branching component of $T$.
\eno
Note that only curves meeting $E$ once are contracted. In particular, $E_0$ is smooth and hence $E_0\simeq \PP^1$.

\bno\label{2.2} Let $t$ denote the number of sprouting contractions in $\psi$. A subdivisional blowing down does not change the quantities $K\cdot(K+D)$ and $E\cdot (K+D)$. Under a sprouting blowing down $K\cdot(K+D)$ increases by 1 and $E\cdot(K+D)$ decreases by 1. Here, by abuse of notation, $K$ denotes the canonical divisor of the image of $\ov S$ at some stage of the contraction process $\psi$, and the images of $E$ and $D$ are denoted by the same letters.  Hence

\med\noin (a) $(E_0+\kn)\cdot (\kn+T)=(E+\ks)\cdot (\ks+D)$.

\med\noin (b) $(E_0+2\kn)\cdot(\kn+T)=(E+2\ks)\cdot(\ks+D)+t=6-2\varepsilon-\ks\cdot E+t=8-2\varepsilon-\gamma+t$.
\eno

We note the following for future reference. \\
\textbf{2.2.1} A contribution (of $1$) to $t$ arises when there is a $(-1)$-curve in $D$ that is non-branching in $D$, meets $E$ once and has attached to it a maximal twig $T$ of $D+E$ consisting of $(-2)$-curves. Note that if $\tau$ their number, then $T$ contributes $\frac{\tau}{\tau +1}$ to $\sum e_i$ in 1.10.\\

\bprop\label{auxiliary}  Suppose that $(E_0+2\kn)(\kn+T)\leq 0$  and that $\ov{N}$ is not
 a Hirzebruch surface or $\PP^2$.  Then there exists a
(-1)-curve $A$ in $\ov{N}$ such that $A\cdot E_0\leq 1$.\eprop

 \begin{proof} Suppose that such a curve does not exist. Let $C_1, C_2$ be the components of $D$ which meet $E$. It may happen that $C_1=C_2$.

\med\noin\textbf{Sub-Lemma 2.3.1.} {\it There is no curve $B$ in
$\ov N$ such that $ (E_0+2\kn)\cdot B<0$}.

\begin{proof} Suppose $B$ exists. Suppose
first  that $|B+\kn+T|\neq \emptyset$. Let $F_m=B+m(\kn+T)$. Arguing as in the proof of Lemma 1.11, we find $m$
 such that $ |F_m|\neq \emptyset$ and $ |F_{m+1}|= \emptyset$ and we have $B+m(\kn+T)=\sum B_i$.
 Then $|B_i+\kn+T|=\emptyset$ for every $i$. By the assumption in Proposition 2.3 we have $0>B\cdot (E_0+2\kn)\geq \sum
 B_i\cdot (E_0+2\kn)$. Hence there exists $B_i$ such that $B_i\cdot (E_0+2\kn)<0$.\\
  Free to replace $B$ by $B_i$, we may assume that $|B+\kn+T|=\emptyset.$
 Then $B$ is a smooth rational curve and $B\cdot T\leq 1$. In particular $B\neq E_0$ since $E_0\cdot T\geq 2$. So $B\cdot E_0 \geq 0$ and $\kn\cdot B<0$,
 i.e., $B^2\geq -1$. Suppose that $B^2\geq 0$. Since $\ov N$ is not a minimal rational surface  there exists a singular member  $\sum B_j$ of $|B|$ such that $B_j^2<0$ for every $j$.
 There exists $B_j$ such that $B_j\cdot (E_0+2\kn)<0$. It follows that $B_j\cdot \kn<0$ hence $B_j^2=-1$, and of course $|B_j+\kn+T|=\emptyset.$ We may replace $B$ by $B_j$. Then $\kn\cdot B=-1$,
 which implies $B\cdot E_0\leq 1$, and $B$ gives a good asymptote for U, a contradiction.
  The sub-lemma is proved.\end{proof}

\med\noin By Theorem 0.1 we have   $\kp (\kn+E_0)=-\infty$. We
argue as in \cite[theorem 2.1]{KM}.

\med (i) Suppose  that $\kn\cdot (\kn+E_0)\leq 0$. Let $L$ be a
(-1)-curve in $\ov N$. Since $L\cdot E_0\geq 2$,  $\vert L+\kn
+E_0\vert\neq \emptyset.$ As above we find $m\geq 1$ such that we have
$$F=L+m(\kn+E_0)=\sum A_i$$
with, for each $i$, $A_i\simeq \PP^1$,
$A_i\cdot E_0\leq 1$ and $A_i^2<0$.
Since $F\cdot\kn<0$, there exists $A_j$ such that $A_j\cdot\kn<0$.
Hence $A_j^2=-1$, so $A_j\cdot(E_0+2\kn)<0$, and we get contradiction with Lemma 2.3.1.

\med (ii) Suppose that $\kn\cdot (\kn+E_0)\geq 1.$ Then
$-\kn-E_0\geq 0$ by the Riemann-Roch theorem and, in fact,
$-\kn-E_0>0$ since $E_0\cdot (-\kn-E_0)=2$. Let again $L$ be a
(-1)-curve in $\ov N$. Write $L=L+\kn+E_0+(-\kn-E_0)$. Since
$h^0(L)=1$ and $L+\kn+E_0\geq 0$, $L+\kn+E_0=0$.
There exists a component $T_i=\psi (D_i)$ of $T$ such that
$T_i\cdot L>0$. Then $T_i\cdot(\kn+E_0)=T_i\cdot (-L)<0$. It
follows that $\kn\cdot T_i<0$. We obtain $T_i\cdot
(E_0+2\kn)=T_i\cdot (E_0+\kn)+T_i\cdot \kn<0$ in contradiction to lemma 2.3.1.

\end{proof}

\bprop\label{basic}  $(E_0+2\kn)(\kn+T)>0$

\eprop

 \begin{proof} We keep notation of 2.1 and 2.2. We have $\kp (\kn+E_0)=-\infty$. Suppose that $(E_0+2\kn)\cdot (\kn+T)\leq 0$. Suppose first   that $\ov N$ is not isomorphic to a  Hirzebruch surface or $\P^2$.
 Let $A$ be a curve as in Proposition 2.3. Suppose that $|A+\kn+T|\neq\emptyset$. We again find $m$ such that $|A+m(\kn+T)\neq\emptyset$ and  $|A+n(\kn+T)=\emptyset$ for $n>m$ and we write
$$F=L+m(\kn+T)=\sum A_i$$
with, for each $i$, $A_i\simeq \PP^1$,
$A_i\cdot T\leq 1$ and $A_i^2<0$.\\
  We have $0>(E_0+2\kn)\cdot (A+m(\kn+T))=\sum (E_0+2\kn)\cdot A_i$. Thus there exists $A_j$ such that $A_j\cdot
 (E_0+2\kn)<0$.  Since $A_0\cdot T\leq 1$, $A_j\neq E_0$. Thus $A_j\cdot
 \kn<0$. Since $A_j^2<0$ we obtain $A_j^2=-1$. It follows that $E_0\cdot A_j\leq 1$. Therefore we may assume that
 $|A+\kn+T|=\emptyset$. Then $A\cdot T\leq 1$, which implies that $A$ is not a branching component of $T$. Also $A\cdot E_0\leq
 1$. By the properties of $T$ it follows that $A$ is not a component of $T+E_0$. The proper transform of $A$ in $\PP^2$ is a good
 asymptote of $U$ and we reach a contradiction.

 \med We have already seen that $\ov N$ cannot be isomorphic to $\PP^2$.
Suppose then that $\ov N$ is isomorphic to a Hirzebruch surface. Since the irreducible components of $T$ generate  $\Pic(\ov
 N)$ freely, $T$ has exactly two components. Write $T=T_1+T_2$. Computing the determinant of $T$ we get $-1=T_1^2T_2^2-1$. We may
 assume therefore that $T_1^2=0$. Let $T_2^2=-n$. Let $a=T_1\cdot E_0, b=T_2\cdot E_0$.
 Then $$E_0\sim (an+b)T_1+aT_2.$$
 Now $p_a(E_0)=0$ implies
 $-2=(an+b)(2a-2)+a(n-2-an)$ and consequently $(a-1)(an+2b-2)=0$. Thus

 \med (i) $a\leq 1$

 or

 (ii) $a\geq 2$ and $2=an+2b.$

 \med In case (i) the proper transform in $\P^2$ of a general member of the system $\vert T_1\vert$ is a good asymptote of  $U$, so (i) cannot occur.\\
  Consider (ii). We have $E_0^2=a^2n+2ab=a(an+2b)=2a.$
  Suppose that $b=0$. Then $2=an$, so $a=2$ and $n=1$. But then $T_2^2=-1$ and $T=T_1+T_2$ is not 2-reduced w.r.t. $E$. Hence $b>0$.  Suppose that $b=1$. Then $an=0$, hence $n=0$, but then the proper transform in $\P^2$ of a general member of the system $\vert T_2\vert$ is a good asymptote of $U.$ Thus $b\geq 2$. Let $q=T_1\cap T_2.$

  \med Suppose that $q\notin E_0$. Then $E_0$ intersects $T_1$ and $T_2$ in points $p_1$ and $p_2$ respectively. The inverse of $\psi$ involves  blowing up over $p_1$ $a$ times and  blowing up over $p_2$ $b$ times. We have $t=2$ (i.e, $\psi $ involves two sprouting contractions w.r.t $D$) and $\kn\cdot(\kn+T)=2+\ks\cdot (\ks+D)=4-\varepsilon -\ks\cdot E$ by 2.2 (b). $\psi$ involves $a+b$ contractions on $E$, hence $\kn\cdot E_0=\ks\cdot E-(a+b)$. We obtain that $\kn\cdot(\kn+T)=4-\varepsilon -\kn\cdot E_0-a-b= 6-\varepsilon +a-b.$ On the other hand $\kn\cdot(\kn+T)=8+\kn\cdot T_1+\kn\cdot T_2=4+n.$ Hence $$n=2-\varepsilon +a-b \leqno (*)$$
  By our assumption
  $$(E_0+2\kn)\cdot (\kn+T)=6-a+b+2n=8-\varepsilon +n\leq 0.
  \leqno (**)$$
  From (ii), since $b=2-\varepsilon+a-n$, we get
  $$(a-2)(n+2)=2\varepsilon-6.
 \leqno(***)$$

$(\star 1)$  By 1.11, $0\leq \varepsilon\leq 3$.

  \med Suppose that $\varepsilon=0$.  $D+E$ has two (-2)-twigs (maximal twigs with each component a $(-2)$ curve) with determinants $a$ and $b$. By \ref{sum ei}, $a=b=2$. From (ii) we get $n=-1$, i.e. $T_2^2=1$. But then the proper transform of $T_2$ in $\ov S$ is a (-1)-curve, so $D+E$ is not NC-minimal w.r.t. $E$.\\

The following four results follow formally from $(\ast)$ (without reference to $q$).\\

$(\star 2)$  $n+2<0$ if $\varepsilon\leq 2$.\\

$(\star 3)$ Suppose that $\varepsilon=1$. We have $n+2=-1$ or -2 or -4, so $n=-3$ or -4 or -6. But $(\ast \ast)$ gives $n\leq -7$.\\

$(\star 4)$ Suppose that $\varepsilon=2$. From $(\ast \ast \ast)$ we get $n+2=-1$ or -2. So $n=-3$ or -4. But $(*)$ gives $n\leq -6.$\\

$(\star 5)$ Suppose that $\varepsilon=3$. From $(\ast \ast)$ we obtain that $n\leq -5$. From $(\ast \ast \ast)$, $(a-2)(n+2)=0$. It follows that $a=2$. \\

Let $T_3$ be the member of the system $\vert T_1\vert$ passing through $p_2$. Since $E_0$ is tangent to $T_2$ at $p_2$, $E_0$ is transversal to $T_3$ at $p_2$. Hence $E_0$ meets $T_3$ transversally at a point $p_3\neq p_2$. After the first blowing up $E_0$ meets the proper transform of $T_2$. It follows that the proper transform of $T_3$ in $\PP^2$ is a good asymptote of $U$.

  \med Now assume that $q\in E_0$.

  \med Assume that $E_0$ meets $T_2$ also in a point $p_2\neq q$. Then $E_0\cap T_1=\{q\}$. Hence $E_0$ is tangent to $T_1$ at $q$. It follows that $E_0$ is transversal to $T_2$ at $q$. Hence the local intersection of $E_0$ with $T_2$ at $p_2$ equals $b-1.$ Suppose that $b=2$. Then $na=-2$, hence $a=2$ and $n=-1$, i.e., $T_2^2=1$. After the first blowing up at $q$, $E_0$ leaves $T_2$. $T_2$ at this stage becomes a 0-curve which meets $E_0$ once. The proper transform of the system $\vert T_2\vert$ in $\P^2$ is a good asymptote of $U$, a contradiction. So $b\geq 3$. It follows that $E_0$ is tangent to $T_2$ at $p_2$. It follows that $\psi $ involves one sprouting contraction w.r.t. $D$. On the other hand, $\psi$ now involves $a+b-1$ contractions on $E$. Computing $\kn\cdot(\kn +E_0)$  as above we get again have $(*)$, hence also $(**)$, $(***)$ and $(\star 1)$ to $(\star 5)$ . \\
  Assume that $\varepsilon=0$. From $(***)$ we get that $n+2$ divides -6. From $(**)$, $n\leq -8$. It follows that $n=-8$, $a=3$. It follows further that $b=13$. The proper transform of $T_1$ in $\ov S$ is a tip of $D+E$ and it is a (-3) curve. $D+E$  also has a twig consisting of 11 (-2)-curves. The twig is created by blowing up over $p_2$. Hence $\sum e_i =\frac{1}{3}+\frac{11}{12}>1$, a contradiction in view of \ref{sum ei}. The cases $\varepsilon=1,2,3$ we eliminate as above.

  \med Assume that $E_0$ meets $T_1$ in a point $p_1\neq q$. Then $E_0\cap T_2=\{q\}.$ Since $b\geq 2$, $E_0$ is tangent to $T_2$ at $q$. Hence $E_0$ is transversal to $T_1$ at $q$. Suppose that $a=2$. Then the proper transform of $T_1$ in $\ov S$ is a (-1)-curve and it meets $E$ once. Thus $D+E$ is not NC-minimal w.r.t. $E$, a contradiction. Hence $a\geq 3$, i.e. $E_0$ is tangent to $T_1$ at $p_1$. As in the previous case $\psi$ involves one sprouting contraction and $a+b-1$ contractions on $E$. Again the $(*)$- and $(\star)$-results hold.\\
Suppose that $\varepsilon=0$. As above we get $n=-8$, $a=3$ and $b=13$. Also $E_0^2=2a=6. $ Let $C_1, C_2$ be the two $(-1)$-components of $D$ which meet $E$. $D-(C_1+C_2)$ has three connected components, two single curves $D_1,D_2$ (they are tips of $D+E$) and one chain $R$. $D_1$ is the proper transform of the curve  produced by the first blowing up over $p_1$ and $D_2$ is the proper transform of $T_2$. $R$ is a chain which has the proper transform of $T_1$ as a tip. It is a (-3)-curve. The rest of $R$ consists of 12 $(-2)$-curves.  We have $D_2^2=-2$, $D_2^2=-5$, $E^2=-9$.\\
 \textbf{2.4.1} Let $Q=D_1+D_2+R+E$. Consider the surface $$Y=\ov S\setminus Q.$$ We claim that $\ovk(Y)\geq 0$. We have $\ks\cdot(\ks+Q)=2+\ks\cdot(\ks+D+E)=4-\varepsilon=4$. Since $Q$ has $4$ components that are rational trees we find $(2\ks+Q)\cdot (\ks+Q) = -\varepsilon =0$.
 By the Riemann-Roch Theorem, $h^0(2\ks+Q)+h^0(-\ks-Q)> 0$. Suppose that $\ovk(Y)=-\infty.$ Suppose  $\vert 2\ks+Q\vert =\emptyset$. Then $-\ks-Q>0$. \\
 \textbf{2.4.2} In view of \ref{kappa S-E}, we have $h^0(-\ks-D-E)=0$ or $-\ks-D-E=0$. Hence by the Riemann-Roch theorem and 1.7
$$h^0(2\ks+D+E)\geq 1+\ks \cdot (\ks+D+E)=3-\epsilon \ \ {\rm or} \ \ \ks=-D-E\ \ {\rm and}\ \ h^0(2\ks+D+E)\geq 2-\epsilon.$$
Hence $2\ks+D+E\geq 0$.  We obtain that $\ks+C_1+C_2=2\ks+D+E+ (-\ks-Q)\geq 0$. This gives  $\ks\geq 0$, a contradiction. Hence $2\ks+Q >0$ and $\ovk(Y)\geq 0$.\\
 \textbf{2.4.3}
  We claim that the pair $(\ov S, Q)$ is almost minimal. If it is not then there exists a (-1)-curve $L$ such that $L\subset\Supp (\ks+Q)^-$  and $L$ is not a component of $Q$. But the intersection matrix of $Q$ is negative definite  and all irreducible components of $Q$ are components of $(\ks+Q)^-$.  Since the rank of $\Pic(\ov S)$ equals the number of irreducible components of $Q$ plus 1 we reach contradiction with \ref{1.13}. \\

  Since   $\chi(Y)=-1$, the BMY-inequality (Langer's version, see 1.13) gives $$\frac{1}{d(D_1)}+\frac{1}{d(D_2)}+\frac{1}{d(R)}+\frac{1}{d(E)}\geq 1.$$
  This is a contradiction since $d(D_1)=2, d(D_2)=5, d(R)=27, d(E)=9$.

  \med The cases $\varepsilon =1,2$ we eliminate as above. If $\varepsilon =3$ we get, as above, that $a=2$, but we we already proved that $a\geq 3.$

  \med Assume that $E_0\cap T=\{q\}$. Then $E_0$ is singular, which we have seen is not the case.
  \end{proof}

 \bcor\label{bound} Let $t$ denote the number of sprouting contractions in $\psi$, see 2.2. Then
 $$7+t\geq 2\varepsilon +\gamma.$$
 \ecor
 \begin{proof}  This follows from 2.2(b) and \ref{basic}.
 \end{proof}

 \section{Separation of branches I: The branches are tangent at infinity }

 \bno Let $\lambda, \tilde\lambda$ be the branches of $\ov U$ at $L_\infty$. The resolution process $\Phi$, see \ref{setup}, can be described in terms of Hamburger-Noether  (HN-) pairs. For the definition in our context  and basic properties of HN-pairs we refer to \cite[1.12]{C-NKR}; see also \cite[Appendix]{KR1} or \cite{Ru1}. We remark also that to each HN-pair there is tacitly associated an $a \in \C$, a parameter that determines the location of the branch on the last exceptional curve produced by the blowups prescribed by the pair. Let $\binom{c_1}{p_1},\dots,\binom{c_h}{p_h}$ (resp.$\binom{\tilde c_1}{\tilde p_1},\dots, \binom{\tilde c_{\tilde h}}{\tilde p_{\tilde h}}$) be the sequence of HN-pairs of $\lambda$ (resp. $\tilde\lambda$). We recall that, by definition, $c_1=\lambda\cdot L_\infty$, $\tilde c_1=\tilde\lambda\cdot L_\infty$, $c_{i+1}=GCD(c_i,p_i)$ and $c_i \geq p_i$.  Let $\mu_1,\mu_2,\cdots$  (resp. $\tilde \mu_1,\tilde
\mu_2,\cdots$) be the sequence of multiplicities of all  singular
points of $\lambda $ infinitely near $\lambda\cap L_\infty$ (resp.  of $\tilde\lambda$ infinitely near  $\tilde\lambda\cap L_\infty$ ).

\med\noin {\bf 3.1.1} Then

$$\sum_{i \geq  1}\mu_i=c_1+p_1+p_2+\cdots +p_h-1.\leqno (i)$$

$$\sum_{i \geq  1}\mu_i^2= c_1p_1+c_2p_2+\cdots +c_hp_h.\leqno (ii)$$

$$\sum_{i \geq  1}\tilde \mu_i=\tilde c_1+\tilde p_1+\tilde
p_2+\cdots +\tilde p_{\tilde h}-1.\leqno(iii)$$

$$\sum_{i \geq  1}\tilde\mu_i^2= \tilde c_1\tilde p_1+\tilde
c_2\tilde p_2+\cdots + \tilde c_{\tilde h}\tilde p_{\tilde h}.\leqno(iv)$$
\eno

Throughout this section we assume that $\lambda\cap L_\infty=\ti\lambda\cap L_\infty=q$ and that the branches cannot be separated by an automorphism of $\C^2$. At the end we will come to a contradiction. We will also assume that the resolution tree $D'$ has the smallest possible number of irreducible components, i.e., if $\sigma\colon \C^2\To\C^2$ is an automorphism, then the number components of the resolution tree of $\ov{\sigma(U)}$ is not less than the number of components of $D'$. \\

\med\noin {\bf 3.1.2} Let $s$ denote the number of {\it common pairs} of $\lambda$ and $\tilde \lambda$. By this we mean that
 $$\frac{c_i}{p_i}=\frac{\tilde c_i}{\tilde p_i} \ \
  {\rm and} \ \ a_i =\tilde a_i \ \ {\rm for}\ \ i=1,\dots, s,$$ but one of these conditions is violated for $i=s+1$. Then the branches separate somewhere along the chains created by the pairs $\binom{c_{s+1}}{p_{s+1}}$, $\binom{\tilde c_{s+1}}{\tilde p_{s+1}}$. Let $m_1,m_2\cdots$ be the sequence of multiplicities of all singular points of $\ov U$ infinitely near to the point $q$.

\med\noin {\bf 3.1.3} We have the following formulas, see [KR], Appendix.

  $$\sum m_i=c_1+\sum p_i-1+\tilde c_1+\sum \tilde p_i-1. \leqno (i)$$

 $$\sum m_i^2=\sum_{i=1}^s (p_i+\tilde p_i)(c_i+\tilde c_i)+\sum_{i>s}p_ic_i+\sum_{i>s}\tilde p_i\tilde c_i+2\min(\tilde p_{s+1}c_{s+1},p_{s+1}\tilde c_{s+1}).\leqno (ii)$$

\ble Let $\gamma'=E'^2$, see 0.1. We obtain the following formulas.

 $$\gamma'+2d=\sum p_i+\sum \tilde p_i.\leqno (a) $$

  $$\gamma'+d^2=\sum_{i=1}^s (p_i+\tilde p_i)(c_i+\tilde c_i)+\sum_{i>s}p_ic_i+\sum_{i>s}\tilde p_i\tilde c_i+2\min(\tilde p_{s+1}c_{s+1},p_{s+1}\tilde c_{s+1}).\leqno (b)$$

\noin where $d=c_1+\tilde c_1.$
\ele
\begin{proof} (a) We have $K_{\P^2}\cdot \ov U=-3d$, $K_{\ov S'}\cdot E'=-2+\gamma'$. Blowing up a point of multiplicity $m$ of  a curve $X$ increases the quantity $K\cdot X$ by $m$. Hence  $K_{\ov S'}\cdot E'-K_{\P^2}\cdot \ov U=-2+\gamma'+3d=\sum m_i$. The statement follows from (i) above.

(b) We have $\ov U^2-E'^2=d^2+\gamma'=\sum m_i^2$. The statement follows now from (ii).
\end{proof}

\med\noin {\bf 3.2.1}  We put  $$c_1-p_1=\alpha c_2, \ \ \tilde c_1-\tilde p_1=\ti\alpha \tilde c_2 \ \ {\rm and} \ \ \alpha_0=\min(\alpha,\ti\alpha).$$

Since $2\min(\tilde p_{s+1}c_{s+1},p_{s+1}\tilde c_{s+1})\leq \tilde p_{s+1}c_{s+1}+p_{s+1}\tilde c_{s+1}$,  3.2(b) gives  $$\gamma'+d^2\leq\sum_{i=1}^s (p_i+\tilde p_i)(c_i+\tilde c_i)+\sum_{i>s}p_ic_i+\sum_{i>s}\tilde p_i\tilde c_i+\tilde p_{s+1}c_{s+1}+p_{s+1}\tilde c_{s+1}$$

\noin and $$\gamma'+d^2\leq (p_1+\tilde p_1)d+(c_2+\tilde c_2)(P+\tilde P),$$

\noin where $\displaystyle P=\sum _{i\geq 2} p_i$, $\displaystyle \tilde P=\sum_{i\geq 2}\tilde p_i.$

\noin From this $$d(d-p_1-\tilde p_1)+\gamma'\leq (c_2+\tilde c_2)(P+\tilde P). \leqno(i)$$

Since $d-p_1-\tilde p_1=c_1-p_1+\tilde c_1-\tilde p_1\geq \alpha_0(c_2+\tilde c_2)$ and since $\gamma'\geq 1$ we get $$\alpha_0 d<P+\tilde P.\leqno(ii)$$

\ble\label{liczba sp}  Let $h_\Phi$ (resp. $h_\Psi$) be the number of sprouting contractions   in $\Phi$ (resp. $\Psi$). Then $h_\Phi=6-\ks\cdot (\ks+D)=2+\varepsilon+\gamma+h_\Psi.$ If $E'$ is not touched by the contractions in $\Psi$, then $\gamma = \gamma'$. If, moreover, $L'_\infty,$ the proper transform of $L_{\infty}$ in $\ov S'$, is not a $(-1)$-curve, then $h_\Psi=0$.
\begin{proof} Under a subdivisional blowing up of a point on a divisor $T$ the quantity $K\cdot (K+T)$ doesn't change. Under a sprouting blowing up the quantity decreases by 1. Hence $h_\Phi= K_{\P^2}\cdot (K_{\P^2}+L_\infty)-K_{\ov S'}\cdot (K_{\ov S'}+D')$ and $h_\Psi =K_{\ov S}\cdot (K_{\ov S} +D) -K_{\ov S'}\cdot (K_{\ov S'}+D')$. Since $K_{\ov S}\cdot E= -2 +\gamma$ the result follows from 1.7.
\end{proof}
\ele

\ble\label{s=0} We have $s=0$.
\ele

\begin{proof}  Suppose that $s\geq 1$. Note that then $\lambda$ and $\tilde \lambda$ are both tangent to $L_{\infty}$. (Otherwise both are not tangent to $L_{\infty}$ and $q$ is a point of multiplicity $deg(\ov U)$ on $\ov U$.) Hence $c_1>p_1, \tilde c_1>\tilde p_1$. Also, both branches have more than one characteristic pair, i.e., $h>1$ and $\tilde h >1$.
We put
$$\frac{c_1}{c_2}= k=\frac{\tilde c_1}{\tilde c_2},\ \ {\rm and} \ \ \frac{p_1}{c_2}=l =\frac{\tilde p_1}{\tilde c_2}.$$
We have $\alpha = \tilde {\alpha} =k-l \geq 1.$

 Suppose that $\alpha=1$, i.e., $k=l+1$. The blowing up over $q$ according to the pair $\binom{l+1}{l}$ produces a chain $L+C+M$, where $L$ has $l$ components with $L_{\infty}$ as a $(-1)$-tip, $C$ is the last exceptional curve and $M$, a $(-l-1)$-curve, is (the proper transform of) the first exceptional curve. The branches $\lambda, \tilde \lambda$ have common center $q'$ on $C \setminus (L\cup M)$. In $\Phi ^{-1}$ we now blow up $q'$. Let $A$ be the resulting exceptional curve. Let us perform $l-1$ successive additional sprouting blowups (they will not be part of $\Phi ^{-1}$), starting with a point on $A$ that is
 {\textit{not}} the center of $\lambda$ or $\tilde \lambda$, creating a chain $A+B$ attached to $C$, with $B$ of length $l-1$. Let $L_{\infty}^\dag$ be the last exceptional curve. As it is well known, we can now blow down, beginning with $L_{\infty}$, the curves in $L$, then $C$, then $A+(B-L_{\infty}^\dag)$, then $M$, producing a new completion of $\C^2$ with $L_{\infty}^\dag$ as new line at infinity and a new completion $\ov U^\dag$ of $U$.\\

\noin\textbf{3.4.1} Let us note for further reference that we have performed an \textit{elementary transformation} of $\C^2$ determined by $q\in L_{\infty}$,  the pair $\binom{l+1}{l}$, the choice of $q'\in C$ and the choice of further fundamental points in creating the chain $B$.\\

 The task of producing the NC-resolution for $\ov U+L_{\infty}$ is accomplished by further blowups over $A$. Let $A^\sharp$ be the resulting configuration of curves and put
 $$\Gamma =L+M+C+A^\sharp+B.$$
Then, as a set, $D'=L+C+M+A^\sharp \subset \Gamma$. In constructing the minimal NC-resolution of $\ov U^\dag+L_{\infty}^\dag$ we have to reconstruct $A$, hence also $B$, and then $A^\sharp$. Hence also $D^\dag \subset \Gamma$. There are three possibilities.\\
(i) The centers of $\lambda$ and $\tilde \lambda$ on $A$ are not on $C$ and $l=1$.
Then $D^\dag=A^\sharp$. (We have $L_{\infty}^\dag=A$.)\\
(ii) The centers of $\lambda$ and $\tilde \lambda$ on $A$ are not on $C$ and $l >1$.
Then $D^\dag=B+A^\sharp + M$.\\
(iii) The center of $\lambda$ or $\tilde \lambda$ on $A$ is on $C$.
Then $D^\dag=B+A^\sharp + M+C$.\\
In each case $D^\dag$ has fewer components than $D'$, contrary to our assumption.\\

 Hence  $\alpha\geq 2$. Then either $L_{\infty}'^2 \leq -2$ or $D'$ has a twig with an initial chain $L_{\infty}'+L$ with $L_{\infty}'^2 =-1$, $L$ a $(-2)$-chain attached to a $(\leq -3)$-curve in $D'$ and $E\cdot (L_{\infty}'+L) = 0$.
 This implies that $E'$ is not touched by the contractions in $\Psi\colon\ov S'\To\ov S$. Hence $\gamma'=\gamma$.\\

 By \ref{bound}, $\gamma\leq 9.$ Suppose that $\alpha\geq 3$. We have  $3d<P+\tilde P$ by 3.2.1(ii). From 3.2(a) we obtain
  $P+\tilde P+ p_1+\tilde p_1-\gamma=2d$. We get that $ d+p_1+\tilde p_1<\gamma\leq 9$. Thus  $d<9-p_1-\tilde p_1\leq 7$. But $d=c_1+\tilde c_1\geq 2(\alpha+1)\geq 8$, a contradiction.

 Hence $\alpha=2$, i.e., $k=l+2.$
 Since $GCD(k,l) =1$, $l$ and $k$ are odd. We have $d-p_1-\tilde p_1=2(c_2+\tilde c_2)$. Substitute this into 3.2.1(i). We obtain $$2d(c_2+\tilde c_2)+\gamma\leq (c_2+\tilde c_2)(P+\tilde P)$$ and $$(c_2+\tilde c_2)(2d-P-\tilde P)+\gamma\leq 0. \leqno(*)$$
 By 3.2(a), $2d-P-\tilde P=p_1+\tilde p_1-\gamma$. Hence $$(c_2+\tilde c_2)(l(c_2+\tilde c_2)-\gamma)+\gamma\leq 0.$$
 From this $l(c_2+\tilde c_2)<\gamma$. Thus $l(c_2+\tilde c_2)\leq 8$, which implies $l\leq 4$. Hence $l \leq 3$. \\
 Suppose that $l=3$. Then $c_2+\tilde c_2=2$ and we have $2(6-\gamma)+\gamma\leq 0$ which gives $\gamma\geq 12$, a contradiction.\\
  Therefore  $l=1$ and $k=3$. $(*)$ takes the form $(c_2+\tilde c_2)(c_2+\tilde c_2-\gamma)+\gamma\leq 0$. By simple algebra we get $(c_2+\tilde c_2-1)(c_2+\tilde c_2+1-\gamma)\leq -1$. This implies that
 $$c_2+\tilde c_2=p_1+\tilde p_1\leq \gamma-2.\leqno(**)$$

 The proper transform of $L_{\infty}$ in $\ov S$ is a $(-2)$-curve. Hence $D'=D$. Blowing up on $L_{\infty}$ according to the pair $\binom{c_1}{p_1}$ produces a chain $L'_{\infty}+C+M$, where  $M$ consists of two $(-2)$-curves and $C$ is branching in $D$.  $L'_{\infty}$ and $M$ are maximal twigs in $D+E$ and contribute $\frac{1}{2}+\frac{2}{3}>1$ to $\sum e_i$. In view of \ref{sum ei} and \ref{bound}, $$\varepsilon\geq 1 \ \ {\rm and}\ \ \gamma \leq 7+t-2\epsilon \leq 5+t.\leqno(***).$$
 Hence $\gamma\leq 7$. If $\gamma=7$, then $\epsilon=1$ and  $t=2$, so $D+E$ has at least two maximal twigs with (-2)-tips and not contained in $L$, those  produced by the pairs $\binom {c_h}{p_h}$ and $\binom{\tilde c_{\ti h}}{\ti p_{\ti h}}$. In view of 1.3.1 we get that  $\sum e_i>2$, in contradiction to \ref{sum ei}.

 Hence we have $\gamma\leq 6$, so $c_2+\tilde c_2\leq 4$.

 Suppose that $c_2=\tilde c_2=2$. Let $a$ be the number of common pairs of type $\binom{2}{2}$. Thus $s\geq a+1$. If  $s\geq a+2$  then the next common pair is of type $\binom{2}{1}$, followed by a number, possibly zero, of common pairs of type $\binom{1}{1}$. The branches then both meet $D$ transversally at different points of the last $(-1)$-curve . It follows that $t=0$ and and $(***)$ gives $\gamma\leq 5$. We reach a contradiction with $(**)$.

  Hence $s=a+1.$ Then either $\binom{c_{s+1}}{p_{s+1}}=\binom{\tilde c_{s+1}}{\tilde p_{s+1}}$, so either $\binom{2}{2}$ or $\binom{2}{1}$, or, say, $\binom{c_{s+1}}{p_{s+1}}=\binom{2}{1}, \binom{\tilde c_{s+1}}{\tilde p_{s+1}}=\binom{2}{2}$. Hence $m=\min(c_{s+1}\tilde p_{s+1}, \tilde c_{s+1}p_{s+1})=2$ or $4$. We have also $d=c_1+\tilde c_1=3(c_2+\tilde c_2)=12$. Note that $\binom{c_{h}}{p_{h}}=\binom{\tilde c_{\tilde h}}{\tilde p_{\tilde h}}=\binom{2}{1}$.
 Hence in 3.2(b) we have $c_hp_h+\tilde c_{\tilde h}\tilde p_{\tilde h}=2+2$ and all other individual terms on the RHS and the term $d^2$ on the LHS  are divisible by $4$. Hence $\gamma$ is divisible by $4$, so $\gamma=4$, and we have a contradiction with $(**)$.

 Suppose now that $c_2=1$. Then $\tilde c_2\leq 3$. Let $a$ be the number of common pairs of type $\binom{\tilde c_2}{\tilde c_2}$. Then $s=1+a$. We have $h=s+1$ with $p_h=c_h=1$ and $\tilde h=1+a+b+1$, where $a+b$ is the total number of pairs equal to $\binom{\tilde c_2}{\tilde c_2}$. We have $m=\min(c_{s+1}\tilde p_{s+1}, \tilde c_{s+1}p_{s+1})= \tilde p_{s+1}$. The formulas 3.2 take the form $$ \gamma +2(3+3\tilde c_2)=s+1+(s+b)\tilde c_2+\tilde p_{\tilde h}. \leqno (1)$$
 and $$\gamma +(3+3\tilde c_2)^2=(3+3\tilde c_2)(1+\tilde c_2)+(s-1)(1+\tilde c_2)^2+2\tilde p_{s+1}+1+b\tilde c_2^2+\tilde c_2\tilde p_{\ti h}.\leqno (2)$$
 Suppose that $\tilde c_2=3$. Then $\gamma=6$ by $(**)$. By $(***)$,  $t\geq 1$. This implies that $\tilde c_{s+1}=\tilde p_{s+1}=3$ and $\tilde p_{\tilde h}=1$. (1) and (2) now give $ 6+24= s+1+3(s+b)+1$, i.e., $28=4s+3b$ and $6+144=48+(s-1)16+9b+10$, i.e., $108=16s+9b$. The system of equations has no integer solutions.\\

Before proceeding with the analysis of cases we note the following. Since $\kappa(\ks+E)=-\infty$ by \cite{Kor}, we have $h^0(2\ks+E)=0$ and by the Riemann-Roch theorem
$$h^0(-\ks -E)\geq \ks \cdot (\ks+E)=\ks^2 -2+\gamma.\leqno(\sharp)$$
As argued in 2.4.2
$$h^0(2\ks+D+E)\geq 1+\ks \cdot (\ks+D+E)=3-\epsilon \ \ {\rm or} \ \ \ks=-D-E. \leqno(\sharp \sharp)$$

 Now suppose that $\tilde c_2=2$. (1) and (2) give  $\gamma+18=s+1 +(s+b)2+1$, i.e., $\gamma+16=3s+2b$ and $\gamma+81=27+(s-1)9+4b+3+2\tilde p_{s+1}  $, i.e., $\gamma+60=9s+4b+2\tilde p_{s+1}$. We have $5\leq\gamma\leq 6$ and $\tilde p_{s+1}\leq 2$. We find two solutions:

 \med (i) $\gamma=5, s=7, b=0, \tilde p_{s+1}=1$

 (ii) $\gamma=6, s=6, b=2, \tilde p_{s+1}=2.$

 In case (i) we find $\ks^2=10-b_2(\ov S)=-2$ and $\varepsilon = 1$. By $(\sharp)$ and $(\sharp\sharp)$, $-\ks-E\geq 0$ and $2\ks+D+E\geq 0$. We obtain $\ks+D=2\ks+D+E+(-\ks-E)\geq 0$, a contradiction. In case (ii) we have $b_2(\ov S)=13$, so $\ks^2=-3$, and $\varepsilon=1$. We come to a contradiction by the same argument.\\

 Suppose that $\tilde c_2=1$. Then $d=6$. The formulas give $\gamma+10=2s+b$ and $\gamma+24=4s+b$. We get the solution

 \med (iii) $\gamma=4, s=7$.

 \med We find $b_2(\ov S)=11$, $\ks^2=-1$, $\varepsilon=2$. Clearly $\ks\neq -D-E$, e.g., $L'_{\infty}\cdot \ks=0,\ \ L'_{\infty}\cdot (-D-E)=1$. We come to a contradiction by the same argument.\\
 \end{proof}

 \bno We have shown that $s=0$. Suppose that the branches stay together after the first blowing up, i.e., they  both are tangent to $L_\infty.$ Let, as in 3.2.1, $$c_1-p_1=\alpha c_2, \ \ \ti c_1-\ti p_1=\ti \alpha \ti c_2.$$
  We will show that, possibly at the cost of increasing the number of components of $D'$, this case can be reduced to the case $\alpha=\ti\alpha=1$ and $\frac{c_1}{p_1}=\frac{\ti c_1}{\ti p_1}=\frac{l+1}{l}$. This case will be dealt with in 3.6.
\eno

 Suppose  that $\ti \alpha=1$. Let $\ti p_1=l\ti c_2$. Then $\ti c_1=(l+1)\ti c_2$. We use the notation of 3.4.1. After blowing up according to $\binom{l+1}{l}$, the center $\ti p$ of $\ti \lambda$ is on $C\setminus (L\cup M)$. We now have three possibilities.\\
 (i) The center $p$ of $\lambda$ is on $M$. Equivalently, $p_1>l(c_1-p_1)$.\\
 (ii)The center $p$ of $\lambda$ is on $C\setminus (L\cup M)$. Equivalently, $p_1=l(c_1-p_1)$, or $\alpha =1$. Moreover, $p\neq p'$.\\
 (iii) The center $p$ of $\lambda$ is on $L$. Equivalently, $p_1<l(c_1-p_1)$.\\

 Suppose we have (i) or (ii). We then perform an elementary transformation exactly as in 3.4.1 with $q'=\ti p$. The argument in the proof of 3.4, slightly modified, shows that we obtain a completion with smaller $D'$. \\

\textbf{3.5.1} We remark that if the HN-sequence for $\tilde{\lambda}$ has at least $l$ pairs $\binom{\ti c_2}{\ti c_2}$ following $\binom{\ti c_1}{\ti p_1}$, or if $\ti c_2=1$, we can construct the above elementary transformation with blowups following $\ti \lambda$, and it will then separates the branches. Hence this is not the case.\\

 Suppose we have (iii). We now perform an elementary transformation as above, but with $q'\neq \ti p$. Then we are in situation (ii) w.r.t. the new coordinate system, i.e., we have  $\alpha=\ti \alpha=1.$ \\

  We may assume that $\alpha\geq 2,\ \ \ti\alpha\geq 2$. We write 3.2(b) as $$\gamma +\alpha c_1c_2+\ti \alpha\ti c_1\ti c_2 +2c_1\ti c_1=\sum\limits_{i\geq 2}p_ic_i+\sum\limits_{i\geq 2}\ti p_i\ti c_i +2\min (c_1\ti p_1, \ti c_1p_1).$$
We have $2c_1\ti c_1=\ti c_1(p_1+\alpha c_2)+c_1(\ti p_1+\ti \alpha \ti c_2)=\ti c_1 p_1+c_1\ti p_1+\alpha c_2\ti c_1+\ti\alpha c_1\ti c_2$.

\noin Therefore $$\gamma +\ti c_1 p_1+c_1\ti p_1+\alpha c_2\ti c_1+\ti\alpha c_1\ti c_2+\alpha c_1c_2+\ti \alpha\ti c_1\ti c_2 =\sum\limits_{i\geq 2}p_ic_i+\sum\limits_{i\geq 2}\ti p_i\ti c_i +2\min (c_1\ti p_1, \ti c_1p_1)$$
and
$$\gamma+\ti c_1 p_1+c_1\ti p_1-2\min (c_1\ti p_1, \ti c_1p_1)+ (c_1+\ti c_1)(\alpha c_2+\ti \alpha\ti c_2) =\sum\limits_{i\geq 2}p_ic_i+\sum\limits_{i\geq 2}\ti p_i\ti c_i .$$

Let $\beta=\ti c_1 p_1+c_1\ti p_1-2\min (c_1\ti p_1, \ti c_1p_1)\geq 0$. We get $$\gamma +\beta +(c_1+\ti c_1)(\alpha c_2+\ti \alpha\ti c_2)=\sum\limits_{i\geq 2}p_ic_i+\sum\limits_{i\geq 2}\ti p_i\ti c_i. \leqno (*)$$
>From 3.2(a) we get $$\gamma+ c_1+\ti c_1+\alpha c_2+\ti\alpha \ti c_2=\sum\limits_{i\geq 2}p_i+\sum\limits_{i\geq 2}\ti p_i. \leqno (**)$$
We may assume by symmetry that $c_2\geq \ti c_2$. Multiply (**) by $c_2$ and subtract (*). We obtain
$$\gamma c_2+(c_1+\ti c_1)c_2 +c_2(\alpha c_2+\ti\alpha\ti c_2)\geq \gamma+\beta +(c_1+\ti c_1)(\alpha c_2+\ti\alpha\ti c_2).$$
So $$\gamma(c_2-1)+(c_1+\ti c_1)c_2\geq \beta +(\alpha c_2+\ti\alpha\ti c_2)(c_1+\ti c_1-c_2).$$
Since $\alpha\geq 2, \ti\alpha\geq 2, \beta\geq 0$ we have

$$\gamma (c_2-1)+(c_1+\ti c_1)c_2\geq ( 2c_2+2)(c_1+\ti c_1-c_2).$$
>From this
$$\gamma (c_2-1)\geq ( 2c_2+2)(c_1+\ti c_1)-( 2c_2+2)c_2-(c_1+\ti c_1)c_2.$$
So
$$\gamma (c_2-1)\geq ( c_2+2)(c_1+\ti c_1)-2( c_2+1)c_2.$$
We have $\gamma\leq 9$ by \ref{bound}, $\ti c_1\geq 3$ since $\ti\alpha\geq 2$. Also $c_1\geq 3c_2$ since $\alpha\geq 2$. We obtain
$$9c_2-9\geq (c_2+2)(3c_2+3)-2(c_2+1)c_2.$$
We get $$0\geq c_2^2-2c_2+15.$$
This is a contradiction.

 \bno In this section we temporarily drop the assumption that $D'$ has the smallest possible number of components. We consider here the case $s=0$, $c_1=(l+1)c_2, p_1=lc_2$, $\ti c_1=(l+1)\tilde c_2, \ti p_1=l\ti c_2$. We will prove that this case does not occur. Suppose opposite.
  Let $H'$ denotes the (-1)-curve produced by the pair $\binom{c_1}{p_1}$ and let $H=\Psi(H')$. The branches meet $H'$ in two different points. $\Psi$ involves $l$ successive contractions beginning with $L_\infty'$.  $H'$ is not contracted by $\Psi$. Let $F$ (resp. $\ti F$) denotes the part of $D'$ produced by the pairs  $\binom{c_2}{p_2},\dots, \binom{c_h}{p_h}$ (resp. $\binom{\ti c_2}{\ti p_2},\dots, \binom{\ti c_{\ti h}}{\ti p_{\ti h}}$). Let $C$ (resp. $\ti C$) be the unique (-1)-curve in $F$ (resp. $\ti F$)
 \eno
  Let $r$ (resp. $\ti r$) denotes the number of pairs equal to $\binom{c_2}{c_2}$ (resp. $\binom{\ti c_2}{\ti c_2}$). Hence $p_{r+2}<c_{r+2}=c_2$ and $c_i\leq {\frac{1}{2}}c_2$ for $i>r+2$. We put $P'=\ds_{i\geq r+2}p_i$. In similar way we define $\ti P'$. Notice that $c_2>1, \ti c_2>1$ by the argument in 3.5.1. Therefore $h>r+1$, $\ti h>\ti r+1$, i.e., $P'\geq 1$ and $\ti P'\geq 1$. Again by 3.5.1 we have $r\leq l-1, \ti r\leq l-1.$

\textbf{3.6.1} We note that $D+E$ has at least 3 maximal twigs, the $-(l+1)$-curve $M$ (see 3.4.1), and one each in $F$ an $\ti F$ with a $\geq (-c_2)$- and a $\geq \ti (-c_2)$-curve as tip respectively. By 1.2.1 they contribute at least $e=\frac{1}{u+1}+\frac{1}{c_2}+\frac{1}{\ti c_2}$ to $\sum{e_i}$ in 1.10. In  particular, $\varepsilon >0$ if $e>1$.

 \med\noin From  3.2(b) we get $$\gamma+d^2=(c_1+\ti c_1)(p_1+\ti p_1)+\ds_{i\geq 2}p_ic_i+\ds_{i\geq 2}\ti p_i\ti c_i=\\
 d(p_1+\ti p_1)+rc_2^2+p_{r+2}c_2+\ti r{\ti c_2}^2+\ti p_{\ti r+2}\ti c_2+\ds_{i\geq  r+3}p_ic_i+\ds_{i\geq \ti r+3}\tilde p_i\ti c_i.$$
 From this $$\gamma+d(c_2+\tilde c_2)\leq rc_2^2+\ti r{\ti c_2}^2+c_2p_{r+2}+\frac{1}{2}c_2(P'-p_{r+2})+\ti c_2\ti p_{r+2}+\frac{1}{2}\ti c_2(\ti P'-\ti p_{\ti r+2}).$$
 From 3.2(a) we get $$d=\frac{p_1+\ti p_1+rc_2+\ti r\ti c_2+P'+\ti P'-\gamma}{2}.$$
  Hence $$\gamma+\frac{1}{2}(c_2+\tilde c_2)(p_1+\ti p_1+rc_2+\ti r\ti c_2+P'+\ti P')-\frac{1}{2}\gamma(c_2+\ti c_2)\leq rc_2^2+\ti r{\ti c_2}^2+\frac{1}{2}c_2p_{r+2}+\frac{1}{2}c_2P'+\frac{1}{2}\ti c_2\ti p_{\ti r+2}+\frac{1}{2}\ti c_2\ti P'.$$

  From this $$(c_2+\ti c_2)(p_1+\ti p_1+rc_2+\ti r\ti c_2)+c_2\ti P'+ \ti c_2P'-\gamma(c_2+\ti c_2)< 2rc_2^2+2\ti r\ti c_2^2+c_2p_{r+2}+\ti c_2\ti p_{\ti r+2}.\leqno (*)$$

  Since $p_1=lc_2, \tilde p_1=l\ti c_2$ and since $P'\geq 1$, $\ti P'\geq 1$ we have
  $$(c_2+\ti c_2)(l(c_2+\ti c_2)+rc_2+\ti r\ti c_2)< 2rc_2^2+2\ti r\ti c_2^2+c_2p_{r+2}+\ti c_2\ti p_{r+2}+(\gamma-1)(c_2+\ti c_2),$$

  $$l(c_2+\ti c_2)^2+(r+\ti r)c_2\ti c_2< rc_2^2+\ti r\ti c_2^2+c_2p_{r+2}+\ti c_2\ti p_{r+2}+(\gamma-1)(c_2+\ti c_2).$$
  Since $r,\ti r\leq l-1$ and $p_{r+2}\leq c_2-1, \ti p_{r+2}\leq \ti c_2-1$

  $$l(c_2+\ti c_2)^2+(r+\ti r)c_2\ti c_2< l(c_2^2+\ti c_2^2)+(\gamma-2)(c_2+\ti c_2).$$
  Finally

  $$c_2\ti c_2(2l+r+\ti r)<(\gamma-2)(c_2+\ti c_2).\leqno (**)$$

  Suppose that $l\geq 3$. Then $6c_2\ti c_2<7(c_2+\ti c_2)$. This implies $c_2=\ti c_2=2$. But then $\varepsilon>0$ by 3.6.1. This implies $\gamma\leq 7$ by \ref{bound}. Now $(**)$ gives $24<20$, a contradiction.\\

  Suppose that $l=2$. Notice that $\gamma<9$. Otherwise $\varepsilon =0$ and $t=2$, and there are two $(-2)$-tips in $D$. This gives a contradiction by \ref{sum ei} as before. $(**)$  gives $4c_2\ti c_2<6(c_2+\ti c_2)$. Let $c_2\leq \ti c_2$.\\
  Suppose that $c_2=2$. We obtain that $8\ti c_2<6(2+\ti c_2)$, i.e., $\ti c_2<6$.  It follows by 3.7.1 that $\varepsilon>0$, so $\gamma\leq 7$. Now $(**)$ gives $8\ti c_2<5(2+\tilde c_2)$, i.e., $\ti c_2\leq 3$. If $\ti c_2=2$, then $\gamma$ is even by 3.2(a), so $\gamma\leq 6$ and $(**)$ gives a contradiction. So $\ti c_2=3$. From $(**)$ we obtain $r=\ti r=0$. We  have $P'=1, \ti P'=\ti p_{\ti r+2}$. Now $(*)$ gives $16<\ti P'$, a contradiction since $\ti P'=1$ or 2.\\
  Suppose that $c_2\geq 3$. Since $\gamma\leq 8$, $(**)$ gives  $c_2(4\ti c_2-6)<6\ti c_2$ and $3(4\ti c_2-6)<6\ti c_2$. We get $\ti c_2<3$, a contradiction.\\

Suppose  $l=1$. Then by 3.5.1 $r=\ti r=0$, i.e $c_2>c_3$ and $\ti c_2>\ti c_3$. We have $d=2c_2+2\ti c_2$ and  the formulas 3.2 take the form

  $$\gamma +3c_2+3\ti c_2=\sum\limits_{i\geq 2}p_i+\sum\limits_{i\geq 2}\ti p_i, \leqno (1)$$

  $$\gamma+2c_2^2+2\ti c_2^2+4c_2\ti c_2=\sum\limits_{i\geq 2}p_ic_i+\sum\limits_{i\geq 2}\ti p_i\ti c_i. \leqno(2)$$
  We may assume that $c_2\geq \ti c_2$. The branches meet the $(-1)$-curve $T_1$ created by $\binom{c_1}{p_1}$ in distinct points. Hence, see 3.3, $h_\Phi = 1+(h-1)+\ti h-1)=h+\ti h-1$. Also, $T_1$ is branching in $D'$ and $(L_\infty')^2=-1$. Hence $\Psi$ contracts only $L'_\infty$, and it is a sprouting contraction, that is $h_\Psi=1$.  By 3.3, $h+\ti h=4+\varepsilon+\gamma$. It follows from \ref{bound} that $\varepsilon+\gamma\leq 8$ ($\gamma=9$ is ruled out as above). Hence $h+\ti h\leq 12$. Since $c_2>1$, $h\geq 2$. Similarly $\ti h\geq 2$. Hence $h,\ti h\leq 10$.

  We write $c_2-p_2=\mu c_3$, $\ti c_2-\ti p_2=\ti\mu \ti c_3$, $c_2=kc_3$, $\ti c_2=\ti k \ti c_3$. Note that $\mu, \ti \mu \geq 1$ and $k, \ti k \geq 2$ since $r, \ti r =0$. We rewrite (2) in the form
  $$\gamma+c_2^2+\ti c_2^2+4c_2\ti c_2=-\mu c_2c_3+\sum\limits_{i\geq 3}p_ic_i-\ti \mu\ti c_2\ti c_3+\sum\limits_{i\geq 3}\ti p_i\ti c_i.\leqno(3)$$
  We get $$\gamma+4c_2\ti c_2\leq c_3^2(h-2-\mu k-k^2)+\ti c_3^2(\ti h-2-\ti \mu\ti k-\ti k^2),$$
  and, since $c_2\geq \ti c_2$,
  $$\gamma\leq c_3^2(h-2-\mu k-k^2)+\ti c_3^2(\ti h-2-\ti \mu\ti k-5\ti k^2).\leqno(4)$$
  We find $\ti h-2-\ti \mu\ti k-5\ti k^2\leq \ti h-24<0$ since $\ti h\leq 10$. It follows from (4) that $$h-2-\mu k-k^2>0.\leqno (5)$$
   Since $h\leq 10$ we get $7\geq k(\mu +k)\geq (\mu+1)(2\mu+1)$. We obtain $\mu=1$ and $k=2$ and $\binom{c_2}{p_2}=c_3\binom{2}{1}$. Hence

\noin  $(***)$ \ \ $D+E$ has at least three tips, two of them $(-2)$-tips. Hence $\varepsilon >0$.\\

( $\binom{c_1}{p_1}$ and $\binom{c_2}{p_2}$ produce $(-2)$-tips, $\binom{\ti c_2}{\ti p_2}$ a third tip.)

  \med\noin\textbf{Claim.} $\gamma+\varepsilon\leq 7$.

  \begin{proof} Suppose otherwise. Then $\varepsilon \geq 2$ is ruled out by 2.5, $\varepsilon =0$ by $(***)$ and 1.10. Hence $\gamma=7, \varepsilon=1$. By 2.5, $t=2$. Suppose that $h>2$. Then $D+E$ has at least four (-2)-tips. It follows from \ref{sum ei} that there are four tips, and they are maximal twigs of $D+E$. Hence $c_h=\ti c_{\ti h}=2$. But now it follows from (2) that $\gamma$ is even, a contradiction. Hence $h=2$. This implies that $c_3=1$, so $c_2=2$. Since $c_2\geq \ti c_2>1$ we have $\ti c_2=2$. We again reach contradiction with (2).\end{proof}

  Since $\gamma+\varepsilon\leq 7$, $h+\ti h\leq 11$. So $h\leq 9$. (5) gives $h>8$. Hence $h=9$ and $\ti h=2$. Also $\gamma+\varepsilon=7$.
  From (2) we get

  $$\gamma+6c_3^2+2\ti c_2^2+4c_2\ti c_2=\sum\limits_{i=3}^9p_ic_i+\ti p_2\ti c_2\leq 6c_3^2+p_9c_9+\ti c_2^2.$$

  It follows that $p_9>1$ since $c_9<4c_2\ti c_2$. We have  $\varepsilon>0$ by $(***)$. Since $\gamma \geq 1$ by 1.8, $\varepsilon \geq 3$ is ruled by 2.5. If $\varepsilon=2$, then $\gamma=5$, so $t=2$ by \ref{bound}, but $p_h =p_9>1$ implies $t\leq 1$.  Hence $\varepsilon=1$ and $\gamma=6$. By \ref{bound} $t\geq 1$. Since $p_9>1$, $\ti p_{\ti h}=\ti p_2=1$. We rewrite (2) and (3) as follows.

  $$5+6c_3+3\ti c_2=\sum\limits_{i=2}^9p_i.\leqno (6)$$

  $$6+8c_3^2+2\ti c_2^2+8c_3\ti c_2=\sum\limits_{i=2}^9p_ic_i+\ti c_2.\leqno(7)$$
  From this
  $$6+6c_3^2+2\ti c_2^2+8c_3\ti c_2=\sum\limits_{i=3}^9p_ic_i+\ti c_2. \leqno(8)$$
  since $c_2=2c_3.$
Suppose that there exists $4\leq j\leq 8$ such that $c_j<c_3$. Then $c_ip_i\leq \frac{c_3^2}{4}$ and $\sum\limits_{i=3}^9p_ic_i \leq (j-3)c_3^2+(10-j)\frac{c_3^2}{4}\leq 6c_3^2$. Now (8) gives a contradiction. Hence $c_i=c_3$ for $i\leq 8$. Suppose that $c_8>p_8$. We write $c_8-p_8=\nu c_9$. Then  $\sum\limits_{i=3}^9p_ic_i\leq 5c_3^2+p_8c_8+p_9c_9=6c_3^2-\nu c_3c_9+p_9c_9\leq 6c_3^2$ and again we reach contradiction with (8). Hence $p_i=c_i$ for $i\leq 8$ and $c_9=c_3$. From (6) we get $$5+3\ti c_2=c_3+p_9.$$
  From (8) we get $$6+2\ti c_2^2+8c_3\ti c_2\leq p_9c_3+\ti c_2.\leqno(9)$$
  Now $p_9=5+3\ti c_2-c_3$ and (9) gives $6+2\ti c_2^2+8c_3\ti c_2\leq (5+3\ti c_2-c_3)c_3+\ti c_2 $. Hence $6+2\ti c_2^2+8c_3\ti c_2+c_3^2\leq 5c_3+3c_3\ti c_2+\ti c_2$, i.e., $$6+2\ti c_2^2+5c_3\ti c_2+c_3^2\leq 5c_3+\ti c_2.$$
  It follows that $6+c_3^2<5c_3$. This gives $2<c_3<3$, a contradiction.

\section{Separation of branches II: The branches separate on the first blowing up}

 In this section we rule out the last case in the proof of theorem 4.16, that of the branches separating on the first blowing up. We assume  that the branch $\lambda$ is tangent to $L_\infty$ and $\ti\lambda$ is not.

 \bno Let $\ti r+1$ denotes the number of pairs of $\ti \lambda $ of the form $\binom{\ti c_1}{\ti c_1}$. So $\ti r\geq 0$. We change slightly our usual labeling. The pairs of $\ti \lambda$ we now label: $$\binom{\ti c_1}{
  \ti c_1},\dots, \binom{\ti c_1}{\ti c_1},  \binom{\ti c_1}{\ti p_1},\dots, \binom{\ti c_{\ti h}}{\ti p_{\ti h}}$$
 with either $\ti c_1=1$ and $\ti r=\ti h=0$, in which case we put $\ti p_1=1$, or $\ti c_1>\ti p_1$. Let $c_1-p_1=\alpha c_2$. We have $\alpha\geq 2$  since otherwise we may, as in 3.6, pass to an embedding with smaller resolution tree.

  Let $T_1$ (resp. $\ti T_1$) be the proper transform in $\ov S$ of the $(-1)$-curve produced by the pair $\binom{c_1}{p_1}$ (resp. $\binom{\ti c_1}{\ti p_1}$). Let $C$ (resp. $\ti C$) be the $(-1)$-curve produced by the last pair in the HN-sequence for $\lambda$ (resp. $\ti \lambda$). Since $\alpha\geq 2$ it is clear that $T_1$, $\tilde T_1$  and $C$, $\ti C$, $E'$ are not touched by $\Psi$, so have the same self-intersection in $\ov S'$ and $\ov S$. In particular $\gamma'=E'^2=E^2=\gamma$.

  Let $S^\ddagger$ be the surface obtained by the first blowup. Let $H^\ddagger, L^\ddagger, E^\ddagger$ be the proper transforms in $S^\ddagger$ of the tangent line $H$ to $\ti\lambda, L_{\infty}, \ov U$. Then $H^\ddagger, L^\ddagger$ are fibers of a $\P^1$-ruling of $S^\ddagger$. We have $E^\ddagger \cdot L^\ddagger =\lambda \cdot L^\ddagger=c_1-p_1$ and $E^\ddagger\cdot H^\ddagger=\ti\lambda\cdot H^\ddagger+ f$, where $f$ is the intersection of $H$ and $\ov U$ at finite distance. We have $f\geq 2$ since otherwise $H$ is a good asymptote.
  If $\ti r=0$, then $\ti\lambda\cdot H^\ddagger=\tilde p_1$. If $\ti r>0$, then $\ti\lambda\cdot H'\geq \ti c_1$. Hence we have the following.
  \eno

  \ble\label{ti r} (a) $c_1-p_1\geq \ti p_1+2\geq 3$.

  \med\noin (b) If $\ti r>0$, then $c_1-p_1\geq \ti c_1+2 \geq 4$.

  \ele

  \bno\label{formulas}  The formulas 3.2 take form
  $$\gamma+2c_1+\ti c_1=\ds_{i\geq 1}p_i+\ti r\ti c_1+\ds_{i\geq 1}\ti p_i\leqno (1)$$
  and $$\gamma+c_1^2+2c_1\ti c_1=\ds_{i\geq 1}p_ic_i+\ti r\ti c_1^2+\ds_{i\geq 1}\ti p_i\ti c_i+2p_1\ti c_1.\leqno (2)$$
  We multiply (1) by $\ti c_1$ and subtract (2). We obtain
  $$\gamma(\ti c_1-1)=(c_1+\ti c_1)(c_1-\ti c_1-p_1)+\sum\limits_{i\geq 2}p_i(\ti c_1-c_i)+\sum\limits_{i\geq 2}\ti p_i(\ti c_1-\ti c_i).\leqno(3)$$
  \eno
  \ble\label{gamma<9} $\gamma\leq 8$.
  \ele
  \begin{proof} Suppose that $\gamma=9$. By \ref{bound}, $\varepsilon=0$ and $t=2$. Hence for both $\lambda$ and $\ti\lambda$ we have the situation described in 2.2.1, i.e., we have two maximal twigs of $D+E$ composed of $(-2)$-curves. If either of these has more than one component, or if $D+E$ has a third maximal twig, we reach a contradiction with \ref{sum ei}. Hence $D+E$ has precisely two maximal twigs, and they are $(-2)$-tips. It follows that $h=\ti h=1$. Let
  $$L_{\infty}--T--T_1$$
  be the upper chain created by the pair $\binom{c_1}{p_1}$, i.e., the chain having $L_{\infty}$ and $T_1$ as tips. Then the chain $L_{\infty}--T$ contracts to a $(-2)$-curve. So either\\
  (i) $L_{\infty}^2=-2$ and $T=\emptyset$ or\\
  (ii) $L_{\infty}^2=-1$  and $T$ has the form $(-2)--\cdots --(-2)--(-3)$ with a number $l\geq 0$ of $(-2)$-curves.\\
  We find  $p_1=1, c_1=3$ in the first case and $p_1=2l+3, c_1=2l+5$ in the second and we reach contradiction with \ref{ti r}(a).
  \end{proof}

  \ble\label{not simple} $\ti c_1>1$, i.e., $\ti\lambda$ is not smooth. In particular, $\ti h \geq 1$ and $\ti c_{\ti h} >\ti p_{\ti h}$.
  \ele
  \begin{proof} Suppose that $\ti c_1=1$. The formulas \ref{formulas}(1) and (2) take the form

   $$\gamma+2c_1+1=p_1+\sum\limits_{i\geq 2} p_i\leqno(1)$$
 and
 $$\gamma+c_1^2+2c_1=p_1c_1+2p_1+\sum\limits_{i\geq 2}p_ic_i.\leqno(2)$$
 We write them in the following form.

 $$\gamma+1+c_1+\alpha c_2=\sum\limits_{i\geq 2}p_i\leqno(3)$$
 and $$\gamma+\alpha c_1c_2+2\alpha c_2=\sum\limits_{i\geq 2}p_ic_i.\leqno(4)$$
 We  multiply (3) by $c_2$ and subtract (4). We get
 $$c_2(1+\gamma) +c_2c_1+\alpha c_2^2\geq \gamma+\alpha c_1c_2+2\alpha c_2.\leqno (5)$$
 From this $$1+\gamma+c_1+\alpha c_2>\alpha c_1+2\alpha.$$
 Let $c_1=kc_2, p_1=lc_2$. Then $\alpha=k-l$.
 We get $$\gamma-2\alpha\geq c_2(k\alpha-\alpha-k).\leqno (6)$$
 Suppose that $\alpha\geq 3$ Then $k=\alpha+l\geq 4$. We obtain $\gamma-6\geq c_2(2k-3)\geq 5c_2$, a contradiction since $\gamma\leq 8.$
 Thus $\alpha=2.$ From (5) we get $$c_2^2(k-2)+c_2(3-\gamma)+\gamma\geq 0.$$
 Therefore $\Delta=(3-\gamma)^2-4\gamma(k-2)\leq 0$. Since $k\geq\alpha+1=3$ we have $(3-\gamma)^2-4\gamma\geq 0$  and finally $\gamma^2-10\gamma+9\geq 0$. From this, since $\gamma>2\alpha=4$ by (6), we obtain $\gamma\geq 9$, a contradiction in view of \ref{gamma<9}.

 \end{proof}

  \ble\label{4.6} $\ti c_1> c_i \ \ {\rm for}  \ \ i\geq 2.$
  \ele
  \begin{proof} It is enough to show that $c_2\geq\ti c_1$ is not possible. Multiply \ref{formulas}(1) by $c_2$ and subtract \ref{formulas}(2). We obtain
  $$\gamma(c_2-1)=-2c_1c_2-\ti c_1c_2+c_1^2+2c_1\ti c_1+p_1c_2-p_1c_1-2p_1\ti c_1+\sum\limits_{i\geq 2}p_i(c_2-c_i)+
  \sum\limits_{i\geq 2}\ti p_i(c_2-\ti c_i)+\ti r\ti c_1(c_2-\ti c_1)+\ti p_1(c_2-\ti c_1).$$
  Let $c_1=kc_2$, $p_1=lc_2$. Then $\alpha=k-l$, hence $k\geq l+2$.

  \med\noin If $c_2\geq\ti c_1$ we get

  $$\gamma(c_2-1)> (-2k+k^2+l-kl)c_2^2+c_2(-\ti c_1+2k\ti c_1-2l\ti c_1).$$
  From this
  $$\gamma>(-2k+k^2+l-kl)c_2-\ti c_1+2k\ti c_1-2l\ti c_1.$$
  Now $-2k+k^2+l-kl=(k-l)(k-2)-l\geq 2(k-2)-l=k+k-l-4\geq k-2\geq 1.$ Since $c_2>\ti c_1$ we obtain
  $$ \gamma>\ti c_1(k^2-kl-l-1)=\ti c_1(k+1)(k-l-1)\geq 4\ti c_1.$$

 Now $\gamma\leq 8$ by \ref{gamma<9}, hence $\ti c_1<2$, a contradiction in view of \ref{not simple}

\end{proof}

\ble\label{beta} Let $\beta=c_1-p_1-\ti c_1$. If $\ti r>0$ then  $2\leq\beta\leq 3$.
\ele
\begin{proof} $\ti r>0$ implies $\beta\geq 2$ by \ref{ti r}(b). By \ref{formulas}(3), \ref{4.6} and 4.2 we find $\gamma(\ti c_1-1)\geq \beta (c_1+\ti c_1) \geq \beta (2\ti c_1 +3)$. In view of 4.6 this gives $\beta < \frac{\gamma}{2}\leq 4$.
\end{proof}

  \bno We consider again the surface $Y$ introduced in 2.4.1.
   Let $Q_1$ (resp. $\ti Q_1$) denote the maximal twig of $D+E$ which meets $C$ (resp. $\ti C$). If $h>1$ then $Q_1$ is the lower subchain produced by the pair $\binom{c_h}{p_h}$. If $h=1$ then $Q_1$ is the image under $\Psi$ of the maximal twig of $D'+E'$ which has  $L'_\infty$ as a tip. In any case $\ti Q_1$ is the lower subchain produced by the pair $\binom{\ti c_{\ti h}}{\ti p_{\ti h}}$. We write $D=Q_1+C+Q_0+\ti C+\ti Q_1$ and put
   $$Q=Q_1+Q_0+\ti Q_1+E \ \ {\rm and} \ \ Y=\ov S\setminus Q.\leqno{\textbf{4.8.1}}$$
   We note $$\chi(Y)= -1.\leqno{\textbf{4.8.2}}$$
  \eno

  \ble\label{k(Y)} If $\gamma\geq 6 $ then $2\ks+Q\geq 0$. In particular $\ovk(Y)\geq 0$.

  \ele

  \begin{proof}   If $\gamma\geq 6$ then $\varepsilon=0 \ \ {\rm or} \ \ 1$ by \ref{bound}. As in 2.4.1 we have $\ks\cdot(\ks+Q)=\ks\cdot(\ks+D+E)-\ks\cdot C-\ks\cdot \ti C=4-\varepsilon.$ If $\varepsilon=0$ we obtain the result as in 2.4.1 and 2.4.2.  \\
   Suppose that $\varepsilon=1$. We have $\ks\cdot(\ks+Q)=3$. By \ref{bound} we have $t\geq 1$. Hence $Q_1$ or $\ti Q_1$, say $\ti Q_1$, consists of (-2)-curves. Then the Riemann-Roch Theorem gives $h^0(-\ks-Q_0-Q_1-E)+h^0(2\ks+Q_0+Q_1+E)> 0$. By 2.4.2 we have $2\ks+D+E\geq 0$. If  $-\ks-Q_0-Q_1-E\geq 0$ then $\ks+C+\ti C+\ti Q_1=2\ks+D+E+(-\ks-Q_0-Q_1-E)\geq 0$. This implies that $\ks\geq 0$, a contradiction. Thus $2\ks+Q_0+Q_1+E\geq 0$ and hence $2\ks+Q\geq 0$.
  \end{proof}

  \ble\label{Y minimal} If $\gamma\geq 6$ then the pair $(\ov S, Q)$ is almost minimal.
  \ele
  \begin{proof}   Suppose that $Q_0$ is contractible (to a quotient singular point), i.e., has negative definite intersection matrix and is a chain or a contractible fork. Then $Q$ has negative definite intersection matrix and the result follows as in 2.4.3.\\
  Suppose that $Q_0$ is not contractible and that $(\ov S, Q)$ is not almost minimal. We need the following.

  \med\noin\textbf{Sublemma} {\it There is no  $(-1)$-curve $L$ in $\ov S$ such that $L\cdot Q_0=0$, $L$ meets two connected components of $Q_1+E+\ti Q_1$ and together with these components contracts to a smooth point.}
  \begin{proof} Suppose that such an $L$ exists. Let $\pi\colon \ov S\To \ov X$ be the contraction of $L$ and the precisely two connected components of $Q_1+E+\ti Q_1$ it meets to a smooth point $q_1$. Let $Q_2$ be the third connected component. The surface $\ov S\setminus Q_0$ is simply connected since it contains $\C^2$. Therefore $X=\ov X\setminus Q_0$ is simply connected. Let $\ov X\To X'$ be the contraction of $Q_2$ to a cyclic singular point $q_2$. Then $X'=\ov X'\setminus Q_0$ is simply connected. It is also easy to compute that $b_2(X')=0$. Hence $X'$ is contractible. Moreover $\ovk(X')=\ovk(X)=\ovk(\ov S\setminus Q_0)=-\infty$. By \cite{KR2} the logarithmic Kodaira dimension of the smooth locus of $X'$ is negative. Since $q_1$ is smooth, $\ovk(X'\setminus \{q_1, q_2\}=\ovk(\ov S\setminus (Q\cup L))-\infty$. It follolws that $\ovk(Y)=\ovk(S\setminus(Q_1+E+\ti Q_1))=-\infty$,  a contradiction in view of \ref{k(Y)}.
  \end{proof}

  Let
$(\ov{Y'},T')$  be an almost minimal model of $(\ov{Y},Q)$. $\ov{Y'}$ is obtained from $\ov{S}$
by a sequence of birational morphisms $p_i : \ol{Y}_i \to \ol{Y}_{i+1}, \ol{S}=\ol{Y}_0 \to
\ol{Y}_1 \to \cdots \to \ol{Y}_\ell=\ol{Y}'$. Let $T_i=(p_{i-1})_*(T_{i-1}), T_0=Q, T'=T_\ell$. Let $Y_i=\ol{Y}_i\setminus T_i$. For every $i$ there exists a $(-1)$-curve $C_i \nsubseteq T_i$ such that
$p_i : \ol{Y}_i \to \ol{Y}_{i+1}$ is the   {\em NC-minimalization} of $C_i+T_i$. Finally, for the almost minimal model $(\ol{Y}',T')$, the negative part
$(K_{\ol{Y}'}+T')^-$ coincides with the bark $\Bk(T')$. The contractions in this process involve only  curves (or their images) contained in the support of
$(K_{\ol{Y}}+T_0)^-$.  We put

$$e(\ol{Y_i},T_i)=\chi(\ol{Y_i}\setminus T_i)+\#\{\mbox{connected components of $T_i$}\}.$$
We find by an elementary calculation that $e(\ov Y_{i+1},T_{i+1})=e(\ov Y_i,T_i)-1.$ Hence $e(\ov Y',T')=e(\ov S,Q)-\ell=3-\ell.$

Let $k$ denote the number of connected components of $T'$ which contract to quotient singularities,
with local fundamental groups $G_j$. Let $u$ denotes the number of connected components of $T'$. By \ref{BMY} we have $$\chi(Y')+\frac{k}{2}\geq \chi(Y')+\sum\limits_{i=1}^k \frac{1}{|G_j|}\geq 0.$$
We have $\chi(Y')=e(\ov Y',T')-u=3-\ell-u$. We obtain $$3-\ell -u +\frac{k}{2}\geq 0.$$
Since $Q_0$ is not contractible, $k\leq u-1$. Also $\ell\geq 1$ since $(\ov S,Q)$ is not almost minimal. We obtain $u\leq 3$ and $k\leq 2$. By the Sublemma above, $\chi(Y')\leq \chi(Y)=-1$ ( in the minimalization process $\chi(Y_{i+1})>\chi(Y_i)$ if and only if $C_i$ meets two connected components of $T_i$ and contracts to a smooth point together with these connected components). From \ref{BMY} we get that $k>1$. Hence $k=2$, $\ell=1$, $u=3$. Also $\chi(Y')=-1=\chi(Y)$. Again by \ref{BMY}\\
 $(*)$ the two contractible connected components of $T'$ are $(-2)$-curves. \\
  We claim that  $C_0$ meets $Q_0$.\\
Suppose otherwise.  Suppose that $C_0$ meets only one connected component $Q_2$ of $Q_1+E+\ti Q_1$. N ow $\chi(Y')=\chi(Y)$ implies that $C_0+Q_2$ must contract to a smooth point. It follows that $Q_2\neq E$ since $E$ is not a $(-2)$-curve (we have $\gamma \geq 6$) and hence $E+C_0$ cannot contract to a smooth point. Therefore $E$ is untouched under $p_0$, so $E^2\neq -2$ in $T'$, and we have a contradiction to $(*)$. Thus $C_0$ meets two connected components of $Q_1+E+\ti Q_1$ and together with these components contracts to a $(-2)$-curve.  Let $\ov X$ be the image of $\ov Y'$ under the contraction of the two connected components of $T'$ that are $(-2)$-curves to singular points. Put $X=\ov X\setminus Q_0$. We have $\ovk(X)=-\infty$, $X$ is simply-connected and has trivial Betti numbers. Hence $X$ is contractible. By \cite{KR2}  the smooth locus of $X$ has negative Kodaira dimension. It follows that $\ovk(Y)=-\infty$, in contradiction to \ref{k(Y)}.\\

Hence $C_0$ meets $Q_0$ and, since $\chi(Y')=\chi(Y)$, one of connected components of $Q_1+E+\ti Q_1. $ The other two connected components are (-2)-curves. It follows that $C_0$ meets $E$ and that $Q_1, \ti Q_1$ are (-2)-curves.

Suppose that $h>1$. Then $d(Q_1)=c_h=2$, $d(\ti Q_1)=\ti c_{\ti h} =2$. By 4.3(2), $4$ divides $\gamma$. Thus $\gamma=8$, $\varepsilon=0$. Now we get contradiction with \ref{sum ei} since $D+E$ has two (-2)-tips and and least one other maximal twig which meets $T_1$.
  Hence $h=1$ and $Q_1$ is a tip of $D+E$ which meets $T_1$.  We reach a contradiction as in the proof of 4.4.

\end{proof}

\bno Put $\omega=h_\Psi$. We have
$h_\Phi =1+\ti r+\ti h+h-1$. By 3.3 we obtain

$$ \ti r+h+\ti h=2+\varepsilon+\gamma+\omega.$$

\ble\label{4.12} If $\gamma\geq 5$ then $Q_0$ is not a chain.
\ele
\begin{proof}

We put $P=\sum\limits_{i\geq 2}p_i$, $\ti P=\sum\limits_{i\geq 2} \ti p_i$. With $\beta$ as in 4.7 we get from \ref{formulas}   $$\gamma(\ti c_1-1)=\beta(c_1+\ti c_1)+\sum\limits_{i\geq 2}p_i(\ti c_1-c_i)+\sum\limits_{i\geq 2}\ti p_i(\ti c_1-\ti c_i).\leqno (*)$$

\med\noin   Suppose that $Q_0$ is a chain. We then have four cases:

\med\noin (a) $h=1$, $\ti h=1$

or

\noin (b) $h=1$, $\ti h=2$, $\ti p_2=1$

or

\noin (c) $h=2$, $p_2=1$, $\ti h=1$

or

\noin (d) $h=2$, $p_2=1$, $\ti h=2$, $\ti p_2=1$.

\med\noin We note the following.\\
(i) If $h=2$, then $\lambda$ produces two tips in $D+E$, one of them a $(-2)$-tip. \\
(ii) If $\ti p_{\ti h}=1$, in particular if $\ti h =2$, then $\ti \lambda$ produces a $(-2)$-tip in $D+E$.

\med We observe that $h+\ti h\leq 4$. From 4.11 we get $\ti r\geq 3$. By \ref{beta} we have $2\leq \beta\leq 3$. Notice that $P+\ti P=h+\ti h-2$. From 4.11 we get $P+\ti P\geq \gamma-\ti r$.
We have $c_1-p_1=\ti c_1+\beta$. From \ref{formulas}(1) we get $$\gamma+c_1+\beta+2\ti c_1\geq \ti p_1+\ti r(\ti c_1-1)+\gamma.$$
So $$c_1+\ti c_1\geq \ti p_1-\beta+\ti r(\ti c_1-1)-\ti c_1.$$
>From $(*)$ we obtain $$\gamma\ti c_1-\gamma\geq \beta(\ti p_1-\beta+\ti r(\ti c_1-1)-\ti c_1)+\sum\limits_{i\geq 2}p_i(\ti c_1-c_i)+\sum\limits_{i\geq 2}\ti p_i(\ti c_1-\ti c_i).\leqno(**)$$

\noin(1) Suppose that $\beta=3$. From 4.11 we have $\ti r\geq\gamma-2$ since $h+\ti h\leq 4$. Using this we get $$2\gamma+3\geq \ti c_1(2\gamma-9)+3\ti p_1.$$
Since  $\beta=3$, $\gamma\geq 7$ by $(*)$. We obtain $17\geq 5\ti c_1+3\ti p_1$. This implies $\ti c_1=2$, $\ti p_1=1$, $\ti h=1$. >From $(**)$ we now obtain $\ti r\leq \gamma-1$. By  4.11, $ 1+h+\gamma-1\geq 2+\varepsilon+\gamma$ hence $h\geq 2+\varepsilon$. Thi gives $\varepsilon=0$, $h=2$. In view of (i) and (ii) we reach contradiction with \ref{sum ei}.\\

\noin (2) Suppose that $\beta=2$.\\
(2.1) Suppose also that $\ti r\geq \gamma-1$. $(**)$ gives $$\gamma+2\geq 2\ti p_1+(\gamma-4)\ti c_1.\leqno{(***)}$$ Since $\gamma\geq 5$ we get $7\geq 2\ti p_1+\ti c_1$. This implies $\ti p_1=1$ or $\ti p_1=2$ and $\ti c_1=3$. In both cases $\ti h=1$.\\
(2.1.1) Suppose also $h=2$ and $\ti p_1=1$. Then $\varepsilon=1$, otherwise we reach contradiction with \ref{sum ei} as above. So $\ti r\geq \gamma$ by 4.11 and $(**)$ gives   $\gamma+4\geq (\gamma-2)\ti c_1+2\ti p_1+p_2(\ti c_1-c_2).$ For $\gamma\geq 6$ we get
 $10\geq 2\ti p_1+4\ti c_1+p_2(\ti c_1-c_2)$, so  $5>\ti p_1+2\ti c_1$ in view of \ref{4.6} and we have  a contradiction since $\ti c_1 \geq 2$. For $\gamma=5$ we get $\ti c_1=2$, $\ti c_1-c_2=1$ and hence $c_2=1$. But then $h=1$.\\
(2.1.2) Suppose also $h=2$ and $\ti p_1=2$. Then $\gamma=5$ by $(***)$. Now $(**)$ gives $15\geq 4\ti r$, but $\ti r\geq \gamma-1=4$, a contradiction.\\
(2.1.3) Suppose also $h=1$. If $\varepsilon+\omega\geq 1$ then 4.11 gives $\ti r\geq \gamma+1$ and $(**)$ gives $\gamma+6\geq 2\ti p_1+\gamma\ti c_1$ and further $11\geq 2\ti p_1+5\ti c_1$; a contradiction.  Hence $\varepsilon=\omega=0$ and $\ti r=\gamma$. >From \ref{formulas}(1) and (3) we get $$\gamma+2c_1+\ti c_1=p_1+\ti p_1+\gamma\ti c_1.$$
and $$ \gamma(\ti c_1-1)=2(c_1+\ti c_1).$$
>From the second equality we have $\gamma\ti c_1=\gamma+2c_1+2\ti c_1$. We substitute it to the first equality and get $$\gamma=p_1+\ti p_1+\ti c_1+\gamma,$$
a contradiction.\\
\noin (2.2) Suppose also that $\ti r\leq\gamma-2$. From 4.11 we obtain $\gamma-2+h+\ti h\geq 2+\varepsilon +\gamma+\varepsilon$, i.e., $h+\ti h\geq 4+\varepsilon+\omega.$ It gives $h=\ti h=2$ and $\varepsilon =0$. We reach contradiction with \ref{sum ei} as before.\end{proof}

\ble\label{not a fork} If $\gamma\geq 6$ then $Q_0$ is not a contractible fork.
\ele
\begin{proof}  Let $H'$  the exceptional curve produced by the first blowing up in $\Phi^{-1}$. Let $H$ denotes the proper transform of $H'$ in $\ov S$.
In view of \ref{not simple} and $c_1>p_1$ we have to blow up at least twice on $H'$. Hence $H^2\leq -3$.

Suppose that $Q_0$ is a fork. Then either $T_1$ or $\ti T_1$ is a branching component in $Q_0$.\\
Suppose $T_1$ is branching. Then the branches are: $R_1$, containing $\Psi(L'_\infty)$;  $R_2$, containing $H$ and $\ti T_1$; $R_3$, meeting $C$. $R_1$ and $R_2$ are maximal twigs of $D+E$.\\
Suppose $\ti T_1$ is branching. Then the branches are: $R_1$ ,containing $\Psi(L'_\infty), T_1, H$;  $R_2$, the lower part of the chain produced by $\binom{\ti c_1}{\ti p_1}$; $R_3$, meeting $\ti C$. $R_1$ and $R_3$ are maximal twigs of $D+E$.\\

Suppose that $Q_0$ is a contractible fork, but not of of type $(2,2,n)$.
Suppose that $T_1$ is a branching component in $Q_0$.

 If $h>2$, then $R_3$ has a $(\leq -3)$-component and hence at most two components. It follows that $h\leq 3$. Also $\ti h=1$ or $\ti h=2$ and $\ti p_2=1$. In any case $h+\ti h\leq 5$. By 4.11, $\ti r\geq 3$. But now the twig $R_2$ has at least 4 components and one of them, $H$, is a $(\leq-3)-$curve. This is impossible.

 Suppose that $\ti T_1$ is a branching component. Again $h+\ti h\leq 5$, so $\ti r\geq 3$. It follows that $R_1$ contains at least 4 components and we reach a contradiction as above.

 Suppose that $Q_0$ is contractible of type $(2,2,n)$.\\ Suppose that $T_1$ is a branching component in $Q_0$. Since $H^2\leq-3$, $R_2$ is the "long" $n$-twig of $Q_0$ and $R_1$, $R_3$ are single (-2)-curves. We have $d(Q_0)=4(n(b-1)-\ti n)$ where $\ti n$ denotes the determinant of the twig $R_2$ with the tip of $R_2$ meeting $T_1$ removed and $b=-T_1^2$. We have $h=2$ and $p_2=2$ since $R_3$ is a single (-2)-curve. So $c_2>p_2$, which implies in particular that $b\geq 3$. Since $R_2$ does not consist of (-2)-curves we have $n-\ti n>1$. We obtain $d(Q_0)\geq 4(2n-\ti n)=4(n+n-\ti n)\geq 4(3+2)=20.$ We have $d(Q_1)=c_2\geq 3$. From \ref{BMY} we get
 $$ 1\leq \frac{1}{d(\ti Q_1)}+\frac{1}{d(Q_1)}+\frac{1}{\gamma}+\frac{1}{d(Q_0)}\leq \frac{1}{d(\ti Q_1)}+\frac{1}{3}+\frac{1}{6}+\frac{1}{20}.\leqno (*)$$
 This implies $d(\ti Q_1)=2$. It follows that $D+E$ has two (-2)-tips. Since it has at least three tips, $\varepsilon=1$ in view of \ref{sum ei}. Thus $\gamma=6 \ \ {\rm or}\ \ 7$. $(*)$ gives $1\leq\frac{1}{2}+\frac{1}{d(Q_1)}+\frac{1}{6}+\frac{1}{20}$, which implies $d(Q_1)\leq 3$. Since $d(Q_1)=c_2\geq 3$ we get $c_2=3$. Moreover , since $\ti c_1>c_2$ by \ref{4.6}, $\ti h=2$.\\

Suppose  that $\omega=0$. Then $R_2=\Psi(L'_\infty)$ and $\frac{c_1}{c_2}=3$, $\frac{p_1}{c_2}=1$. Hence $c_1=9,\ \ p_1=3$.  From 4.11 we obtain $\ti r= \gamma-1.$ Now \ref{formulas}(1) gives $\gamma+12=\ti p_1+(\gamma-2)\ti c_1$. Since $\ti h=2$, $\ti c_2\geq 4$ and $\ti p_1\geq 2$. We obtain $\gamma+12\geq 2+4(\gamma-2)$. This implies $\gamma\leq 6$, so $\gamma=6$. Also $\ti c_1=4$ and $\ti p_1=2$, $\ti r=5$. From \ref{formulas}(3) we obtain $6(4-1)=13\cdot 2+2(\ti c_1-3)+\ti c_1-2$, a contradiction.\\

 Thus $\omega\geq 1$.  From 4.11 we now get $\ti r\geq \gamma$. We have $c_1-p_1=\ti c_1+\beta$, so \ref{formulas}(1) gives $\gamma+c_1+\beta+2\ti c_1\geq \ti p_1+\gamma\ti c_1+2+1$, i.e., $c_1\geq \ti p_1+(\gamma-2)\ti c_1+3-\gamma-\beta$. Now \ref{formulas}(3) gives $$\gamma\ti c_1-\gamma\geq (\ti p_1+(\gamma-1)\ti c_1+3-\gamma-\beta)\beta+\ti c_1-2+2(\ti c_1-3).$$
 If $\beta=3$ then $2\gamma+8\geq 3\ti p_1+2\gamma\ti c_1$. If $\beta =2$, then $\gamma+6\geq 2\ti p_1+(\gamma+1)\ti c_1$. In both cases we get contradiction since $\gamma=6$ or $7$ and $\ti c_1\geq 4.$\\

 Assume that $\ti T_1$ is a branching in $Q_0$. Now $R_2$ and $R_3$ are single (-2)-curves. Hence $\ti h=2$, $\frac{\ti c_1}{\ti c_2}=2$ and $\ti p_2=2$. We again have $d(Q_0)\geq 20$ and, by \ref{BMY}, we get $d(Q_1)=2$ and $d(\ti Q_1)=\ti c_2=3$. Hence $\ti c_1=6, \ti p_1=3$. From \ref{formulas}(3) we get $$5\gamma=(c_1+6)\beta +\ti p_2(\ti c_1-\ti c_2)+p_2(\ti c_1-c_2).$$
 Suppose that $h=1$. Then $5\gamma=(c_1+6)\beta+6$. As above, $\gamma=6$ or $7$. Since $\beta=2$ or $3$, $\beta$ divides $\gamma$. Hence $\gamma=6$. Now $30\geq 2c_1+18$, which gives $c_1\leq 6=\ti c_1$. But $c_1=p+1+\ti c_1+\beta>\ti c_1$. We reach a contradiction.\\ Suppose that $h=2$. Then $c_2=d(Q_1)=2$ and $p_2=1$. We get $5\gamma=(c_1+6)\beta+2(\ti c_1-\ti c_2)+\ti c_1-2=(c_1+6)\beta+10.$  Since $\gamma\leq 7$, we have $35\geq2c_1+12+10$  and again $c_1\leq 6$, a contradiction.

\end{proof}

\bprop $\gamma\leq 5.$
\eprop
\begin{proof} Suppose that $\gamma\geq 6$. By 4.12 and 4.13, $Q_0$ is not contractible. By \ref{BMY} we have $\frac{1}{d(Q_1)}+\frac{1}{d(\ti
Q_1)}+\frac{1}{\gamma}\geq 1$. Since $\gamma\geq 6$ we have

\med\noin (i) $d(Q_1)=d(\ti Q_1)=2$

or

\med\noin (ii) $\{d(Q_1), d(\ti Q_1)\}=\{2,3\}$, $\gamma=6$. In this
case $\ovk(Y)=0$ or $1$ since otherwise $((K_{\ov X}+D)^+)^2>0$  in 1.13. We record that $(BkE)^2 = -\frac{4}{\gamma}$. Put $B_0=(\Bk Q_1)^2+(\Bk \ti Q_1)^2$. Then $B_0=-4$ if $Q,\ \ \ti Q$ consist of $(-2)$-curves. Otherwise they are single curves and $B_0=-\frac{4}{3}-\frac{2}{2}=-\frac{10}{3}$.

\med Consider (i). Then $\ti c_{\ti h}=d(\ti Q_1)=2$.  Suppose that $h>1$.
Then $c_h=d(Q_1)=2$. By \ref{formulas}(2), $4$ divides $\gamma$. Hence
$\gamma=8$. So $\varepsilon=0$ by \ref{basic} and we reach
contradiction with \ref{sum ei}. Suppose that $h=1$. Then $C=T_1$ and $Q_1$ contains $\Psi(L'_\infty)$. We come to contradiction as in the proof of 4.8.

\med Consider (ii). By 4.9, $2\ks+Q\geq 0$. Let
$\ks+Q=P+\Bk Q$ be the Zariski decomposition. We have
$P\cdot(\ks+Q)=P^2=0$ since $\ovk(Y)=0$ or $1$. Recall that $P$ is {\it nef}. We get
$$0=P\cdot(2\ks+2Q)=P\cdot(2\ks+Q)+P\cdot Q\geq P\cdot Q.$$
Hence $P\cdot Q=P\cdot Q_0=0$. Fujita \cite{Fu1} classifies connected
components $Q_0$ of a boundary divisor of an almost minimal surface
such that $P\cdot Q_0=0$. In our case $Q_0$ is one of the following:

\med\noin (a) a chain,

\med\noin (b) a tree with exactly two branching components and four
maximal twigs being $(-2)$-tips,

\med\noin (c) a fork of type $(d_1,d_2,d_3)$ where
$\frac{1}{d_1}+\frac{1}{d_2}+\frac{1}{d_3}=1$.

\med Case (a) is ruled out by \ref{4.12}.\\
\med Consider (b). Then $(\Bk
Q_0)^2=-2$.  Now $-4-\varepsilon=(\ks+Q)^2=(\Bk
Q)^2=-2-\frac{4}{\gamma}+ B_0$ and $\frac{4}{\gamma}+B_0$ is an integer. Since $\gamma\geq 6$ this implies $B_0=-\frac{10}{3}$ and $\gamma=6$. We get $-4-\varepsilon=(\Bk
Q)^2=-2-4=-6$, which gives $\varepsilon=2$, a contradiction by
\ref{bound}.

\med Consider (c).
 We have $$-4-\varepsilon=(\ks+Q)^2=(\Bk Q_0)^2+B _0-\frac{4}{6}.\leqno{(*)}$$
 $Q_0$ is of the type (3,3,3), (2,4,4) or (2,3,6). We find that  $(\Bk Q_0)^2\leq -1$. Since $\varepsilon\leq 1$ it follows from $(*)$ that $\varepsilon=1$. It follows next from $(*)$ that $B_0=-\frac{10}{3}$, i.e, that $Q_1$ and $\ti Q_1$ are single curves, and that $(\Bk Q_0)^2=-1$. By examining all possibilities we see that every twig of $Q_0$ is a tip, i.e., $\# Q_0=4$. Hence $\#Q=7$, $b_2(\ov S)=8$, $\ks^2=2$. From $\ks\cdot (\ks+Q)=5$ we get $\ks\cdot Q=-7$.  Let $B$ be the branching component of $Q_0$. Examining all possibilities we find that $B^2>0$. But $B=T_1$ or $B=\ti T_1$ and both $T_1$ and $\ti T_1$ are untouched under $\Psi$ and hence are negative curves. We reach a contradiction.

\end{proof}

\ble $\gamma-\ti c_1-p_1-\ti p_1>0$.\ele

\begin{proof} Suppose the opposite. By 4.14, $\gamma\leq 5$, so $4\geq \ti c_1+p_1+\ti p_1$. In view of \ref{not simple} we get $\ti c_1=2, \ti p_1=1$ and $p_1=1$. It follows that $\ti h=1$ and $c_2=1$. Hence $h=1$. By 4.11, $\ti r\geq \gamma+\varepsilon +\omega.$ Hence $\ti r>0$, otherwise $\gamma=0$, which is impossible by 1.8. By \ref{ti r}, $2\leq\beta\leq 3.$ \ref{formulas}(3) gives $5\geq \gamma=(c_1+2)\beta$. It follows that $c_1=0$, a contradiction.

\end{proof}
\bthm If $U$ has no good asymptote then the branches of $\ov U$ at infinity can be separated by an automorphism of $\C^2$.
\ethm

\begin{proof} Suppose opposite. By results of section 3 we may assume that things are as in 4.1. By 4.14 we have $\gamma\leq 5. $

 From \ref{formulas}(2) we get
   $$\gamma+\alpha c_2c_1+2\alpha c_2\ti c_1=\ds_{i\geq 2}p_ic_i+\ds_{i\geq 2}\ti p_i\ti c_i+\ti r\ti c_1^2+\tilde p_1\ti c_1.\leqno(1)$$

  Since $\ti p_1\leq \alpha c_2-2$
  $$\gamma+\alpha c_2c_1+\alpha c_2\ti c_1+2\ti c_1\leq\ds_{i\geq 2}p_ic_i+\ds_{i\geq 2}\ti p_i\ti c_i +\ti r\ti c_1^2.\leqno (2)$$

  \noin \ref{formulas}(1) takes the form $$\gamma+2c_1+\ti c_1-p_1-\ti p_1=\ds_{i\geq 2}p_i+\ds_{i\geq 2}\ti p_i+\ti r\ti c_1 \leqno (3)$$ and
  $$2(c_1+\ti c_1)+\gamma-\ti c_1-p_1-\ti p_1=\ds_{i\geq 2}p_i+\ds_{i\geq 2}\ti p_i+\ti r\ti c_1.$$

  By 4.15, $\gamma-\ti c_1-p_1-\ti p_1\leq 0$. Hence $\displaystyle c_1+\ti c_1\geq\frac{1}{2}(\ds_{i\geq 2}p_i+\ds_{i\geq 2}\ti p_i)+\frac{1}{2}\ti r\ti c_1.$
  From (2) we get  $$\gamma+\frac{\alpha c_2}{2}(\ds_{i\geq 2}p_i+\ds_{i\geq 2}\ti p_i+\ti r\ti c_1)+2\ti c_1\leq \ds_{i\geq 2}p_ic_i+\ds_{i\geq 2}\ti p_i\ti c_i +\ti r\ti c_1^2.$$

  Suppose that  $\ti r=0$. Then

  $$\gamma+\frac{\alpha c_2}{2}(\ds_{i\geq 2}p_i+\ds_{i\geq 2}\ti p_i)+2\ti c_1\leq \ds_{i\geq 2}p_ic_i+\ds_{i\geq 2}\ti p_i\ti c_i .$$
  Since $\alpha\geq 2$ this implies that $\frac{\alpha c_2}{2}<\ti c_2$ and further $c_2<\ti c_2$.
  It follows that $\ti c_2=\ti p_1$, otherwise $\ti c_2\leq \frac{\ti p_1}{2}\leq\frac{\alpha c_2-2}{2}.$ We rewrite (3) as
  $$\gamma+c_1+\alpha c_2+\ti c_1-\ti c_2= \ds_{i\geq 2}p_i+\ds_{i\geq 2}\ti p_i$$ and (1)  as $$\gamma+2\alpha c_2\ti c_1+\alpha c_1c_2-\ti c_2\ti c_1=\ds_{i\geq 2}p_ic_i+\ds_{i\geq 2}\ti p_i\ti c_i.$$

  Multiply the first equality by $\ti c_2$ and subtract the second one. We obtain
  $$\gamma(\ti c_2-1)=(\alpha c_2-\ti c_2)( c_1+2\ti c_1-\ti c_2)+\sum_{i\geq 2}p_i(\ti c_i-c_i)+\sum_{i\geq 2}\ti p_i(\ti c_2-\ti c_i).$$

  Since $\ti c_2>\frac{\alpha c_2}{2}\geq c_2\geq c_i$ for $i\geq 2$ $$\gamma(\ti c_2-1)\geq (\alpha c_2-\ti c_2)( c_1+2\ti c_1-\ti c_2)$$

  We have $\ti c_1\geq 2\ti c_2.$ Since $c_1>\ti p_1=\ti c_2$  and $\alpha c_2-\ti c_2=c_1-p_1-\ti p_1\geq 2$ we obtain $\gamma(\ti c_2-1)\geq 2\cdot 4\ti c_2$. It follows that $\gamma=9$, a contradiction.

  Thus $\ti r>0$. By \ref{beta}, $2\leq\beta\leq 3$.

  Suppose that $\ti c_1\geq 2c_2$. Then $\ti c_1-c_i\geq \frac{\ti c_1}{2}$ for every $i\geq 2$. Also $\ti c_1-\ti c_i\geq \frac{\ti c_1}{2}$ for every $i\geq 2$.

  From \ref{formulas}(3)  we obtain

  $$\gamma(\ti c_1-1)=(c_1+\ti c_1)\beta+\sum\limits_{i\geq 2}p_i(\ti c_1-c_i)+\sum\limits_{i\geq 2}\ti p_i(\ti c_1-\ti c_i)\geq (c_1+\ti c_1)\beta+\frac{\ti c_1}{2}(\sum\limits_{i\geq 2}p_i+\sum\limits_{i\geq 2}\ti p_i).\leqno(4)$$
  It follows that $\gamma\geq 5$, i.e., $\gamma=5 $, and further $5>4+\frac{1}{2}(\sum\limits_{i\geq 2}p_i+\sum\limits_{i\geq 2}\ti p_i)$. It gives $\sum\limits_{i\geq 2}p_i+\sum\limits_{i\geq 2}\ti p_i\leq 1.$ It follows that $h=1$ or $h=2$ and $p_2=1$, and similarly for $\ti h,\ \ \ti p_1$. It follows that $Q_0$ is a chain in contradiction to \ref{4.12}.

  Hence $\ti c_1<2c_2.$ (4) and \ref{4.6} give $5\ti c_1> 2(c_1+\ti c_1)$ i.e. $c_1<\frac{3}{2}\ti c_1$. But $c_1\geq 3c_2>\frac{3}{2}\ti c_1$ since $\alpha\geq 2$, a contradiction.

  \end{proof}

\eno

\end{document}